\documentclass[12pt]{article}
\usepackage{mathrsfs}
\usepackage{amsfonts}
\usepackage{graphicx}
\begin{document}
\title{Boundary reducible handle additions on simple 3-manifolds}

\author {Yannan Li, Ruifeng Qiu and Mingxing Zhang}

\date{}
\maketitle
\begin{abstract} Let $M$ be a simple manifold, and $F$ be a
component of $\partial M$ of genus two. For a slope $\gamma$ on $F$,
we denote by $M(\gamma)$ the manifold obtained by attaching a
2-handle to $M$ along a regular neighborhood of $\gamma$ on $F$. In
this paper, we shall prove that there is at most one separating
slope $\gamma$ on $F$ so that $M(\gamma)$ is $\partial$-reducible.
\end{abstract}


\section{Introduction}
\ \ \ \ \ Let $M$ be a compact, orientable 3-manifold such that
$\partial M$ contains no spherical components. $M$ is said to be
simple if $M$ is irreducible, $\partial$-irreducible, anannular and
atoroidal.

Let $M$ be a simple 3-manifold, and $F$ be a component of $\partial
M$. A slope $\gamma$ on $F$ is an isotopy class of essential simple
closed curves on $F$. For a slope $\gamma$ on $F$, we denote by
$M(\gamma)$ the manifold obtained by attaching a 2-handle to $M$
along a regular neighborhood of $\gamma$ on $F$, then capping off a
possible 2-sphere component of the resulting manifold by a 3-ball.

The distance between two slopes $\alpha$ and $\beta$ on $F$, denoted
by $\Delta(\alpha,\beta)$, is the minimal geometric intersection
number among all the curves representing the slopes.

In this paper, we shall study $\partial$-reducible handle additions
on simple 3-manifolds. The main result is the following theorem:

{\bf Theorem 1.} \ Let $M$ be a simple 3-manifold with $F$ a genus
two component of $\partial M$, then there is at most one separating
slope $\gamma$ on $F$ such that $M(\gamma)$ is $\partial$-reducible.

Using the same method, we can prove the following theorem:

{\bf Theorem 2.} \ Suppose that $M$ is a simple 3-manifold with $F$
a genus at least two component of $\partial M$, and  $\gamma_1$ and
$\gamma_2$ are two separating slopes on $F$. If $\partial M-F$ is
compressible in each of $M(\gamma_1)$ and $M(\gamma_2)$, and
$\partial M(\gamma_i)-(\partial M-F)$ are incompressible in
$M(\gamma_i)$ for $i=1,2$, then $\Delta(\alpha,\beta)\leq 2$.

{\bf Comments on Theorem 1 and Theorem 2.} \

1. If $F$ is a torus, then $M(\gamma)$ is the Dehn filling along
$\gamma$. Y. Wu has shown that there are at most three slopes
$\gamma$ on $F$ so that $M(\gamma)$ is $\partial$-reducible. In this
case, $\gamma$ is non-separating on $F$. But if $g(F)>1$, then it is
possible that there are infinitely many non-separating slopes
$\gamma$ on $F$ so that $M(\gamma)$ is $\partial$-reducible.

2. Suppose that $g(F)>1$. Scharlemann and Wu[SW] proved that there
are only finitely many basic degenerating slopes on $F$. As a
corollary of this result, there are only finitely many separating
slopes such that $M(\gamma)$ is not simple. Recently, we([ZQL])
proved that $\Delta(\alpha,\beta)\leq 4$ when $M(\alpha)$ and
$M(\beta)$ are reducible.





\section {Preliminary}

\ \ \ \ \ Let $M$ be a simple 3-manifold with $F$ a genus two
component of $\partial M$. In the following arguments, we assume
that $\alpha$ and $\beta$ are two separating slopes on $F$ such that
$M(\alpha)$ and $M(\beta)$ are $\partial$-reducible. We denote by
$P$ and $Q$ the $\partial$-reducing disks of $M(\alpha)$ and
$M(\beta)$.


{\bf Lemma 2.1[SW].} \ If $\Delta(\alpha,\beta)>0$, then each of
$M(\alpha)$ and $M(\beta)$ is irreducible.\hfill$\Box$\vskip 5mm

{\bf Lemma 2.2.} \ $\partial P$ and $\partial Q$ are disjoint from
$F$.

{\bf Proof.} \ Without loss of generality, we assume that $\partial
P\subset F$. Then $\partial P$ lies in one of the two toral
components of $\partial M(\alpha)$ produced from $F$. This means
that $M(\alpha)$ is reducible, contradicting Lemma 2.1.
\hfill$\Box$\vskip 5mm

{\bf Lemma 2.3.} \ If $\alpha \neq \beta$, then
$\Delta(\alpha,\beta)\geq 4$.

{\bf Proof.} \ This follows the assumptions that $g(F)=2$ and
$\alpha$ and $\beta$ are separating. \hfill$\Box$\vskip 5mm

{\bf Lemma 2.4.} \ There is an incompressible and
$\partial$-incompressible planar surface, say $S_\alpha$(resp.
$S_\beta$) in $M$ with all boundary components but one having the
same slope $\alpha$(resp. $\beta$).


{\bf Proof.} \ By definition, $M(\alpha)=M\cup_{N(\alpha)}D\times
[0,1]$, where $N(\alpha)$ is a regular neighborhood of $\alpha$ on
$F$. Now let $P$ be a $\partial$-reducing disk of $M(\alpha)$. By
Lemma 2.2, $\partial P\subset \partial M-F$. Now each component of
$P\cap (D\times [0,1])$ is a disk $D\times \{t\}$ ($t\in[0,1]$). We
may assume that $|P\cap (D\times [0,1])|$ is minimal among all
$\partial$-reducing disks of $M(\alpha)$. Hence $S_\alpha=P\cap M$
is incompressible in $M$.

Suppose, otherwise, that $S_\alpha$ is $\partial$-compressible with
$B$ a $\partial$-compressing disk. Let $\partial B = u\cup v$, where
$u$ is an essential arc in $S_\alpha$, and $v$ is an arc in
$\partial M$. Since $S_\alpha$ is incompressible, $v$ is essential
on $\partial M-\partial S_\alpha$.

(1) \ $v$ has endpoints on the different components of $\partial
S_\alpha$.

Now $\partial$-compressing $S_\alpha$ along $B$ will give a planar
surface with fewer boundary components. It is a contradiction.

(2) $v$ has endpoints on the same component of $\partial S_\alpha$.

Now either $\partial v\subset \partial P$ or $\partial v\subset C$,
where $C$ is a component of $\partial S_\alpha-\partial P$. If
$\partial v\subset \partial P$, then $\partial$-compressing
$S_\alpha$ along $B$ will give a planar surface with fewer boundary
components. If $\partial v\subset C$, then, by $\partial$-reducing
$S_\alpha$ along $B$, we can obtain a $\partial$-reducing disk of
$M(\alpha)$, say $P^{'}$, such that $\partial P^{'}\subset F$.
Contradicting Lemma 2.2. \hfill$\Box$\vskip 5mm

The components of $\partial S_\alpha$(resp. $\partial S_\beta$)
lying on $F$ are called inner components of $\partial
S_\alpha$(resp. $\partial S_\beta$). We denote by $p$ and $q$ the
numbers of inner components of $\partial S_\alpha$ and $\partial
S_\beta$. Number the inner components of $\partial S_\alpha$(resp.
$\partial S_\beta$) by $\partial_u S_\alpha$(resp. $\partial_i
S_\beta$) for $u=1,2,\cdots,p$(resp. $i=1,2,\cdots,q$), such that
they appear consecutively on $\partial M$. This means that
$\partial_u S_\alpha$ and $\partial_{u+1} S_\alpha$ bound an annulus
in $\partial M$ with its interior disjoint from $S_\alpha$.

Isotopy $S_\alpha$ and $S_\beta$ so that $|S_\alpha \cap S_\beta|$
is minimal.

{\bf Lemma 2.5.} \ Each component of $S_\alpha \cap S_\beta$ is
essential on both $S_\alpha$ and $S_\beta$. \hfill$\Box$\vskip 5mm

{\bf Proof.} \ This follows from lemma 2.4. \hfill$\Box$\vskip 5mm

By Lemma 2.2, $\partial P$(resp. $\partial Q$) is disjoint from the
inner components of $\partial S_\beta$(resp. $\partial S_\alpha$).

Let $\Gamma_P$ be a graph on $P$ obtained by taking the arc
components of $S_\alpha \cap S_\beta$ as edges and taking the
components of $\partial S_\alpha$ as fat vertices. The inner
components of $\partial S_\alpha$ are called inner vertices.
Specially, $\partial P$ is called the outer vertex of $\Gamma_P$.
Similarly, we can define $\Gamma_Q$ on $Q$.

Let $\Gamma\in\{\Gamma_P,\Gamma_Q\}$. An edge $e$ in $\Gamma$ is
called an inner edge if the two endpoints of $e$ are incident to
inner vertices of $\Gamma$, $e$ is called a boundary edge if one of
the endpoints of $e$ is incident to the outer vertex of $\Gamma$.

In this section, the definitions of a cycle, the length of a cycle,
a disk face and parallel edges are standard, see [GL], [SW] and
[Wu].

{\bf Lemma 2.6.} \ There are no 1-sided disk faces on
$\Gamma_P$(resp.$\Gamma_Q$). \hfill$\Box$\vskip 5mm


{\bf Lemma 2.7[SW].} \ There are not common parallel edges in both
$\Gamma_P$ and $\Gamma_Q$.\label{lemmaparallel}

{\bf Proof.} \ The proof follows from Lemma 2.1 of
[SW].\hfill$\Box$\vskip 5mm


\begin{center}
\includegraphics[totalheight=5.5cm]{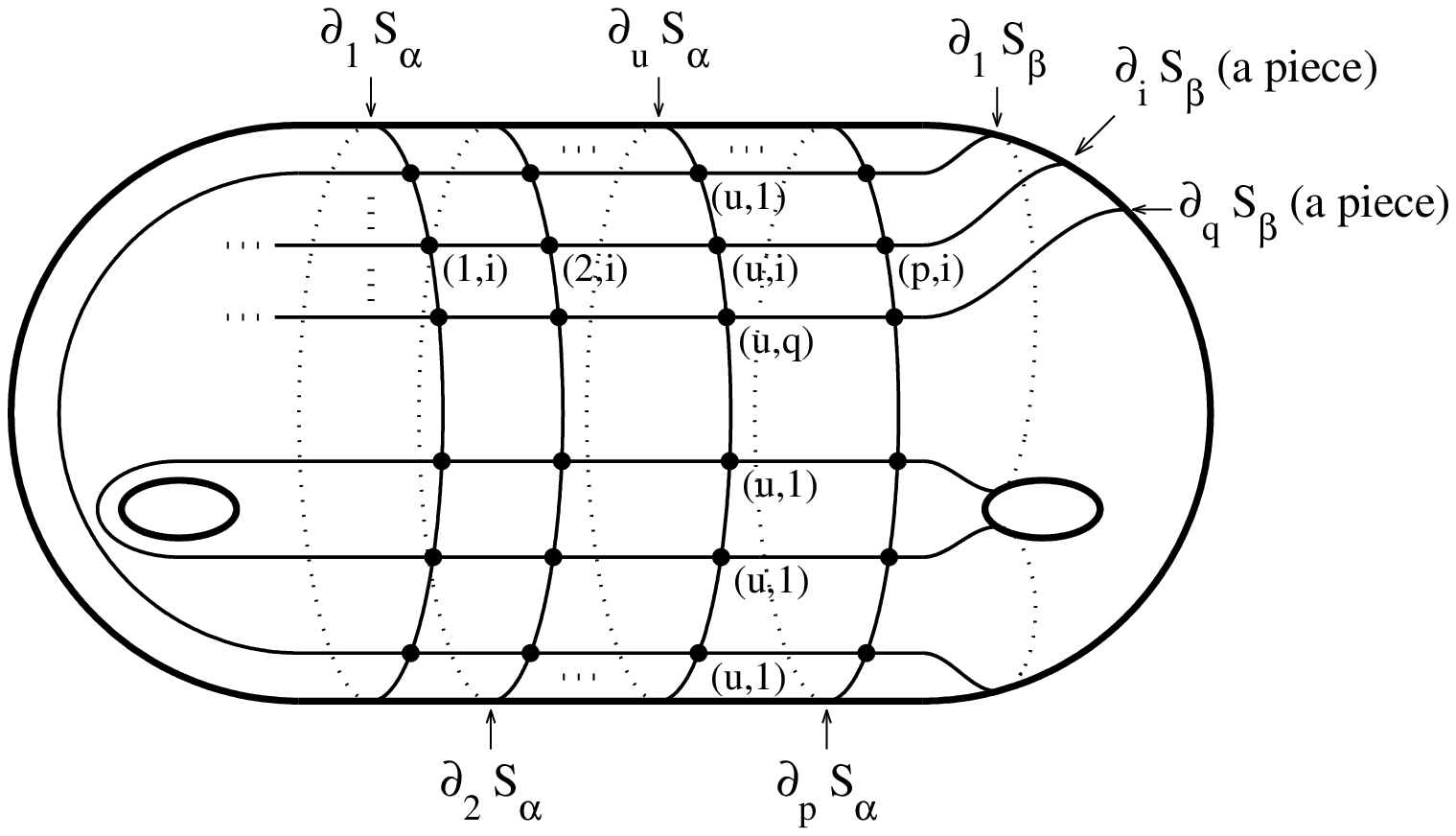}
\begin{center}
Figure 1.
\end{center}
\end{center}

Let $x$ be an endpoint of an edge lying on the inner vertices of
$\Gamma_P$ and $\Gamma_Q$. If $x\in \partial_u
S_\alpha\cap\partial_i S_\beta$, then $x$ is labeled $(u,i)$(see
Figure 1), or $i$(resp. $u$) in $\Gamma_P$(resp. $\Gamma_Q$) for
shortness when $u$(resp. $i$) is specified(see Figure 2). Now when
we travel around $\partial_u S_\alpha$, the labels appear in the
order $1, 2,\cdots,q, q, \cdots, 2, 1,\cdots$(repeated
$\Delta(\alpha,\beta)/2$ times). If $x \in \partial P\cap \partial
Q$, then $x$ is labeled $*$.

Now each edge has a label pair induced by the labels of its two
endpoints. That is to say, each inner edge $e$ of $\Gamma_P$(resp.
$\Gamma_Q$) can be labeled with $(u,i)-(v,j)$, or $i-j$(resp. $u-v$)
in $\Gamma_P$(resp. $\Gamma_Q$) for shortness; a boundary edge can
be labeled with $(u,i)-*$.

\begin{center}
\includegraphics[totalheight=2.5cm]{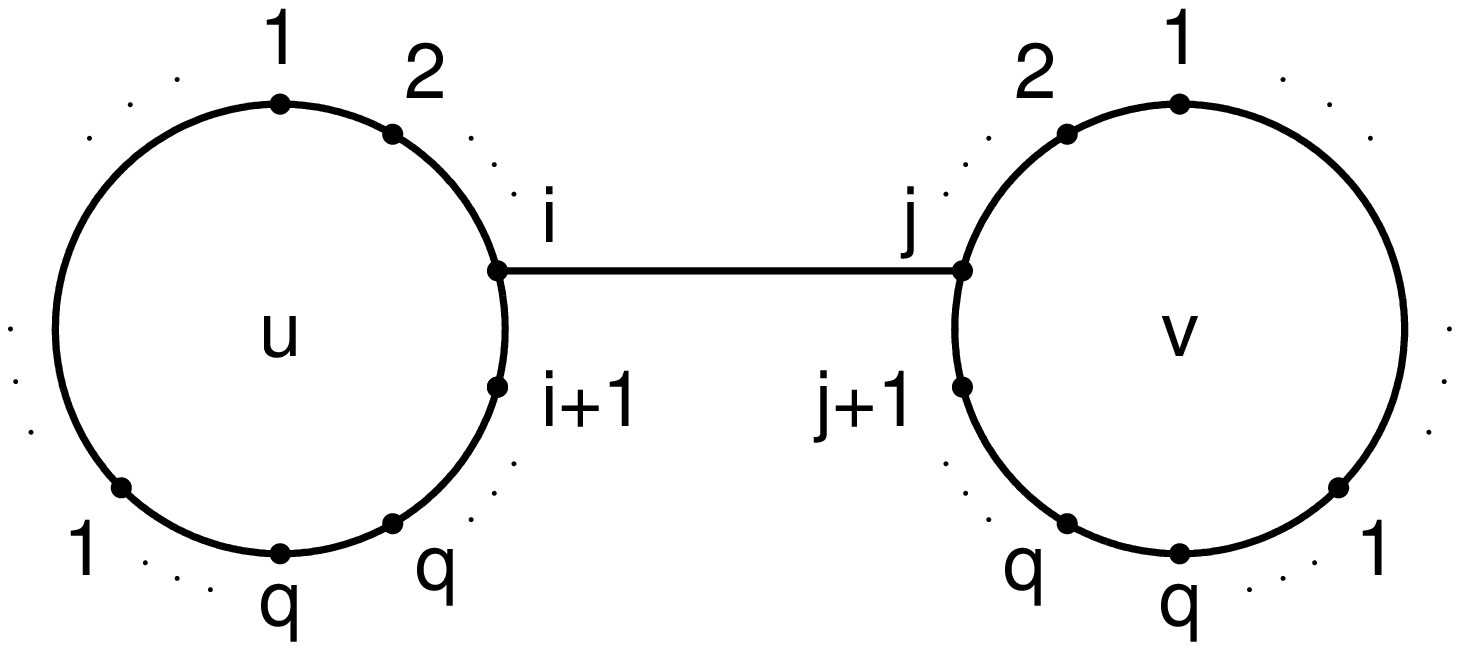}
\includegraphics[totalheight=2.5cm]{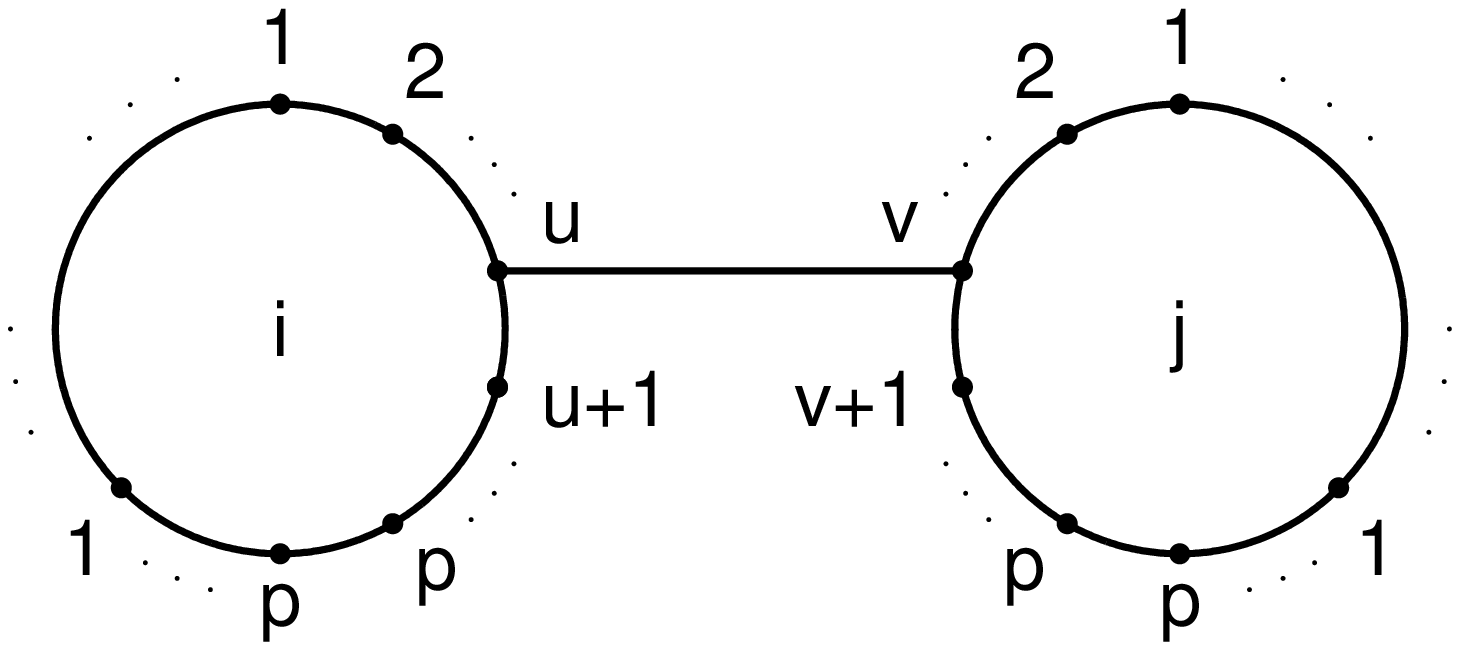}
\begin{center}
Figure 2: Labels on $\Gamma_P$ and $\Gamma_Q$.
\end{center}
\end{center}

Let $C=\{e_k \ | \ k=1,2,\cdots,n\}$ be a cycle in $\Gamma$ such
that $\partial_2 e_{k}=\partial_1 e_{k+1}$ and $\partial_2
e_{n}=\partial_1 e_1$. $C$ is called a virtual Scharlemann cycle if
$C$ bounds a disk face in $\Gamma$, and $e_k$ is labeled $i-j$(or
$u-v$) for each $1 \leq k \leq n$. Furthermore, if $i\neq j$(or
$u\neq v$), then $C$ is called a Scharlemann cycle.

{\bf Lemma 2.8[CGLS].} \ $\Gamma$ contains no Scharlemann
cycle.\hfill$\Box$\vskip 5mm

{\bf Lemma 2.9.} \ $\Gamma_P$(resp. $\Gamma_Q$) contains no
$2q$(resp. $2p$) parallel edges.

{\bf Proof.} \ Suppose, otherwise, $\Gamma_P$ contains $2q$ parallel
edges $e_1,e_2,\cdots,e_{2q}$. Then, by Lemma 5.2 in [ZQL],
$e_1,e_2,\cdots,e_{2q}$ are boundary edges. For each $1 \leq i\leq
q$, there are just two edges, each of which is labeled $i-*$. The
two edges form a length two cycle in $\Gamma_Q$ connecting  the
vertex $\partial_i S_\beta$ to $\partial Q$. Then the two edges in
the innermost one of these cycles are parallel in $\Gamma_Q$,
contradicting Lemma 2.7. \hfill$\Box$\vskip 5mm

\section {Parity rule}
\ \ \ \ \ By Lemma 2.3, we may assume that $\Delta(\alpha,\beta)\geq
4$.

Fix the directions on $\alpha$ and $\beta$. Then each point in
$\alpha \cap \beta$ can be signed ``$+$'' or ``$-$'' depending on
whether the direction determined by right-hand rule from $\alpha$ to
$\beta$ is pointed to the outside of $M$ or to the inside of $M$.
See Figure 3. Since $\alpha$ and $\beta$ are separating, the signs
``$+$'' and ``$-$'' appear alternately on both $\alpha$ and $\beta$.
For details, see [ZQL].

\begin{center}
\includegraphics[totalheight=8cm]{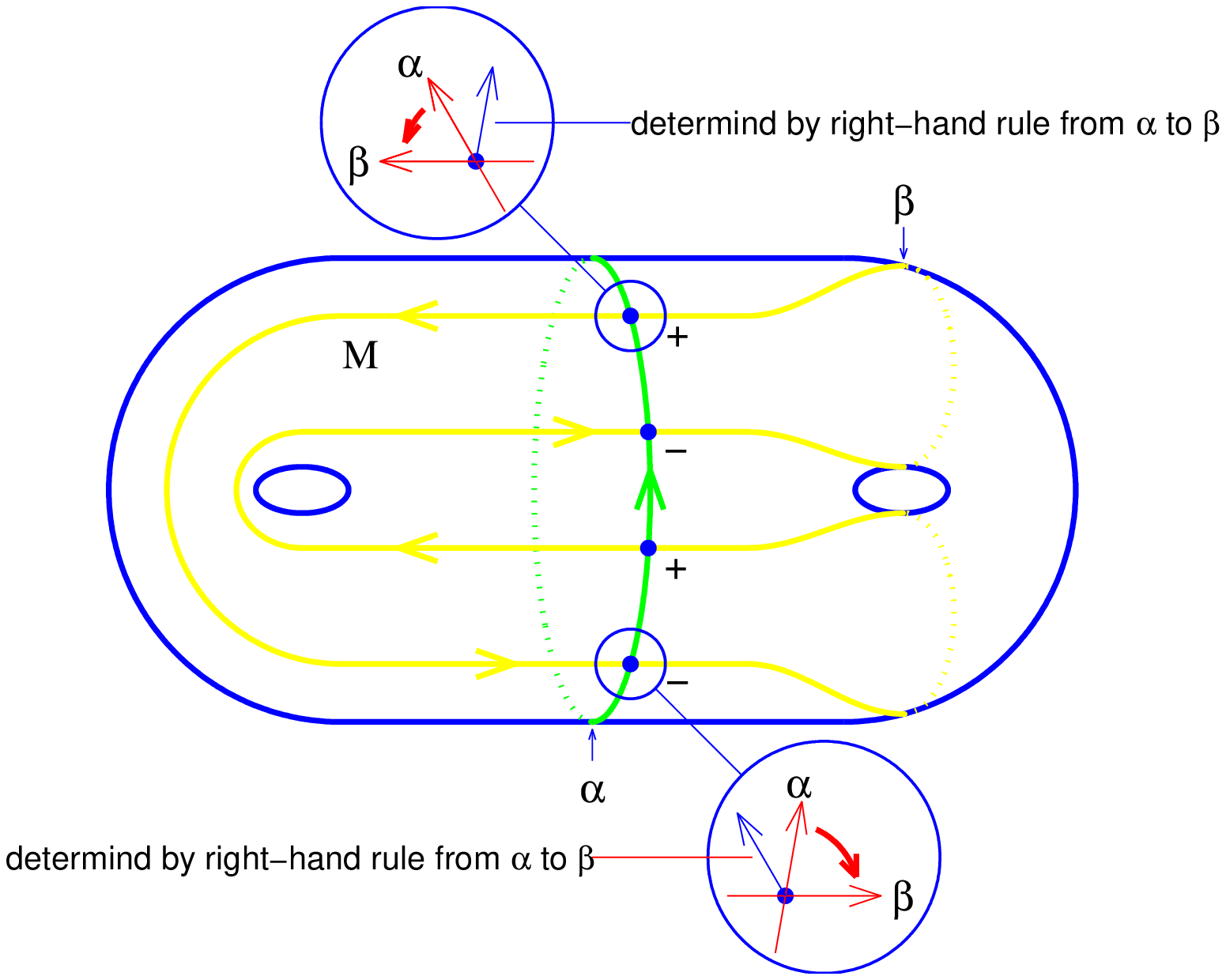}
\begin{center}
Figure 3
\end{center}
\end{center}

\begin{center}
\includegraphics[totalheight=6cm]{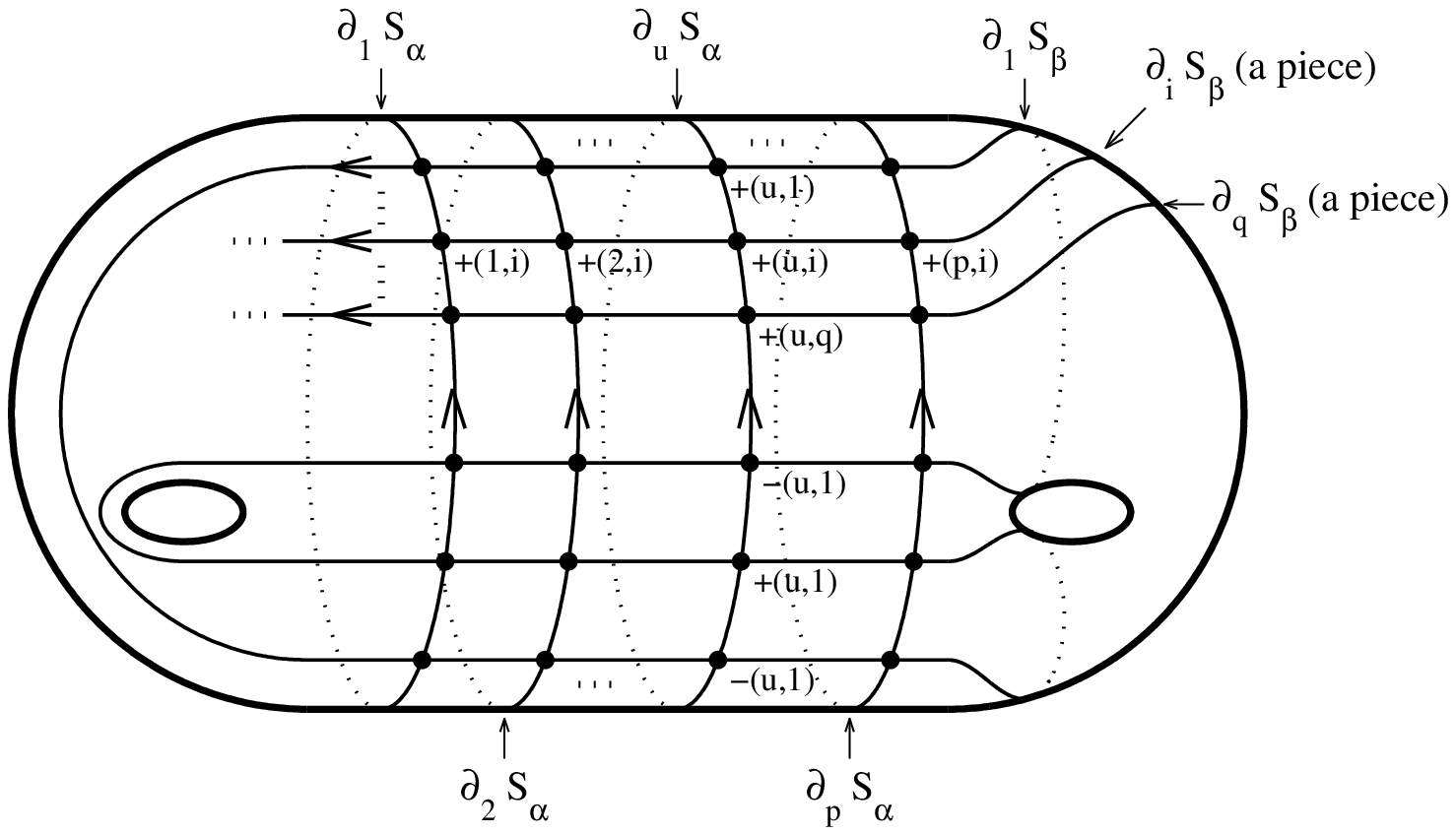}
\begin{center}
Figure 4: Signs on inner endpoints
\end{center}
\end{center}

Give a direction to each inner component of $S_\alpha$(resp.
$S_\beta$) such that they are all parallel to $\alpha$(resp.
$\beta$) on $\partial M$. Then each inner endpoint $x\in
(S_\alpha-\partial P)\cap(S_\beta-\partial Q)$ can be signed as
above. We denoted by $c(x)$ the sign of $x$. See Figure 4. Now the
signed labels appear on $\partial_u S_\alpha$ as
$+1,+2,\cdots,+q,-q,\cdots,-2,-1,\cdots$, (repeated
$\Delta(\alpha,\beta)/2$ times). See Figure 5.

\begin{center}

\includegraphics[totalheight=4 truecm]{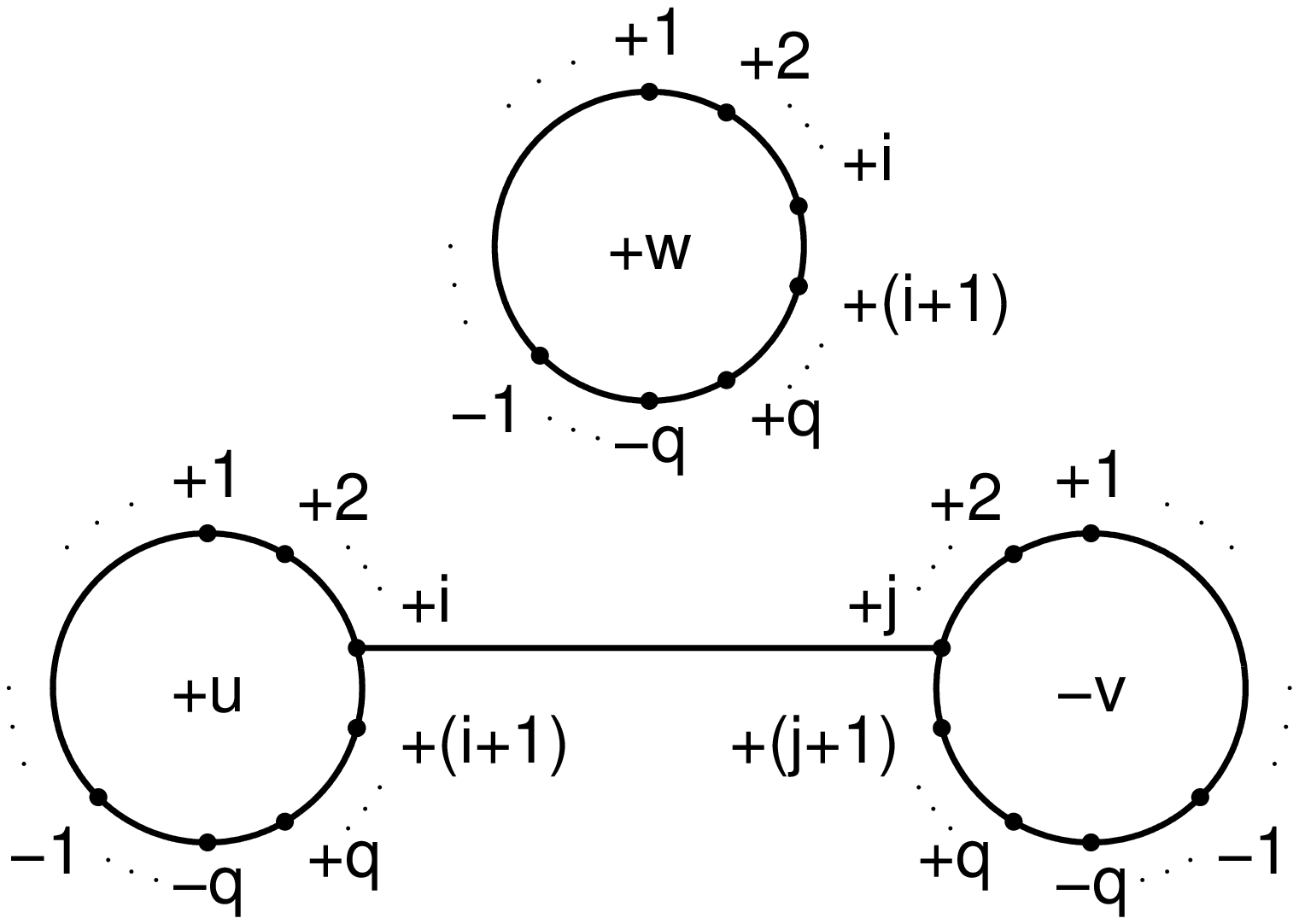}~~~~
\includegraphics[totalheight=4 truecm]{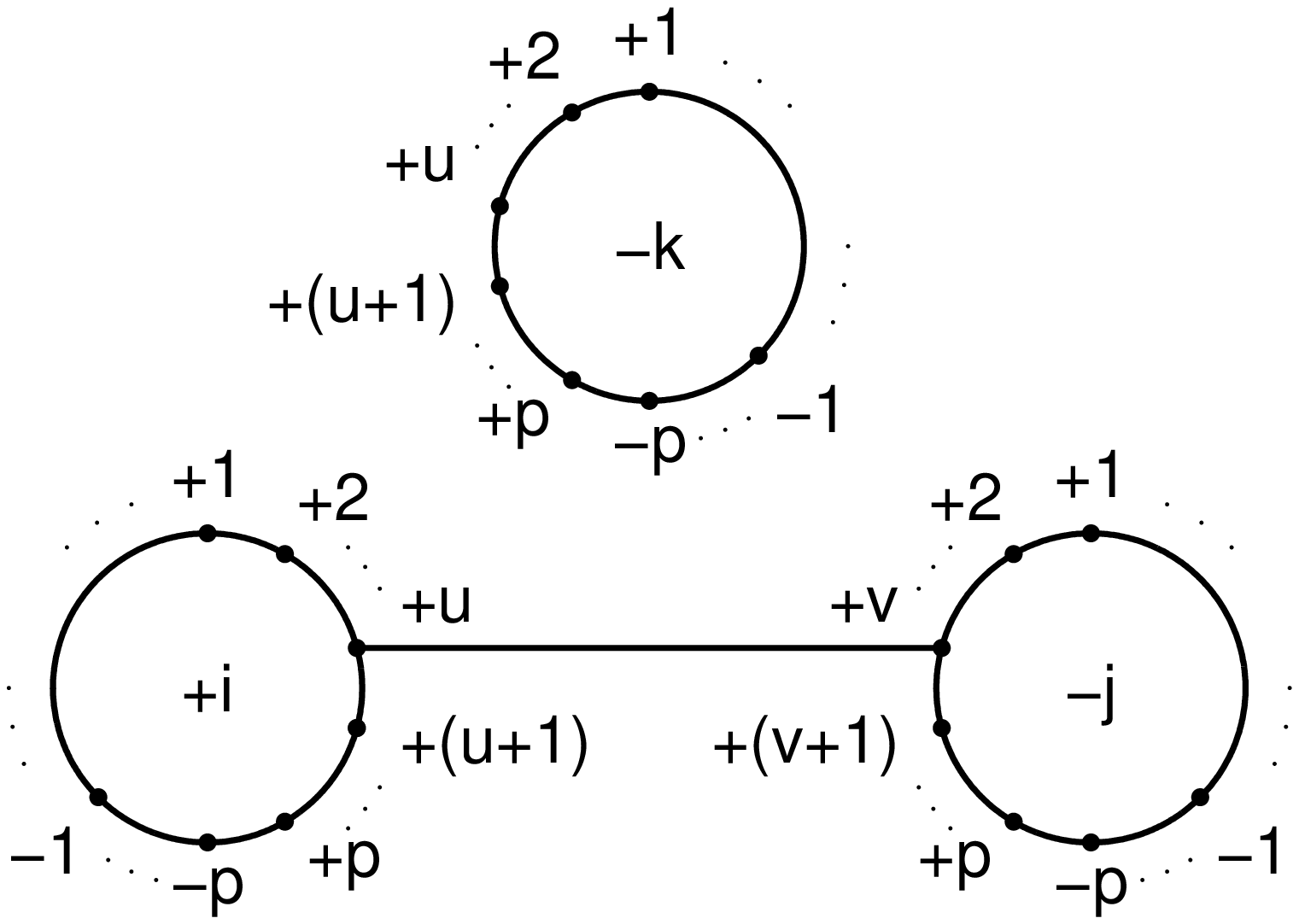}

labels on $\Gamma_P$~~~~~~~~~~~~~~~~labels on $\Gamma_Q$
\begin{center}
Figure 5
\end{center}
\end{center}

Now we sign the inner vertices of $\Gamma_P$. Suppose $P\times
[0,1]$ is a thin regular neighborhood of $P$ in $M$. Let
$P^+=P\times 1$ and $P^-=P\times 0$. For some $1\leq u\leq p, 1\leq
i\leq q$, let $c$ be a component of $\partial_u S_\alpha\times[0,1]
\cap
\partial_i S_\beta$ with the induced direction of $\partial_i S_\beta$. We
define $s(u)$ the sign of $\partial_u S_\alpha$ as follows:

(1) Suppose $c$ intersects $\partial_u S_\alpha$ at a ``$+$'' point,
we define $s(u)=+$(resp. $s(u)=-$) if the direction of $c$ is from
$P^+$ to $P^-$(resp. from $P^-$ to $P^+$).

(2) Suppose $c$ intersects $\partial_u S_\alpha$ at a ``$-$'' point,
we define $s(u)=+$(resp. $s(u)=-$) if the direction of $c$ is from
$P^-$ to $P^+$(resp. from $P^+$ to $P^-$).

Since for each $1\leq i\leq q$, $\partial_i S_\beta$ has the same
direction with $\beta$ on $F$, the definition as above is
independent of the choices of $c$ and $i$.

Since $M$ is orientable, $\partial_u S_\alpha$ and $\partial_v
S_\alpha$ have the same direction on $P$ when $\partial_u S_\alpha$
and $\partial_v S_\alpha$ have the same signs. This means the labels
$+1,+2,\cdots,+q$, $-q,\cdots,-1$ of the inner endpoints appear on
both $\partial_u S_\alpha$ and $\partial_v S_\alpha$ are in the same
direction in $\Gamma_P$. Similarly, the labels
$+1,+2,\cdots,+q,-q\cdots,-1$ appear in opposite the directions when
$\partial_u S_\alpha$ and $\partial_v S_\alpha$ have different
signs. See Figure 5.

And we sign the inner vertices of $\Gamma_Q$ in the same way as
$\Gamma_P$.

The labels with the signs defined as above are said to be type A.
Now we have Parity rule A:

{\bf Lemma 3.1[ZQL].} \ For an edge $e$ in $\Gamma_P$(and
$\Gamma_Q$) with its endpoints $x$ labeled $(u,i)$ and $y$ labeled
$(v,j)$, the following equality holds:

\hfill $s(i)s(j)s(u)s(v)c(x)c(y)=-1$ ~~~~~ ($*$).\hfill$\Box$

On $\Gamma_P$ and $\Gamma_Q$, we define new signs of inner endpoints
of edges as follows:

For each inner endpoint $x$ in $\Gamma_P$(resp. $\Gamma_Q$) with
labele $(u,i)$. Let $g(x)=c(x)\times s(u)$(resp.
$g^\prime(x)=c(x)\times s(i)$). Then the signed label $g(x)i$ (resp.
$g^\prime(x)u$) on $\Gamma_P$(resp. $\Gamma_Q$) of $x$ is said to be
type B.

{\bf Remark 1.} Under type B labels, the labels $+1$, $+2\cdots$,
$+q,-q,\cdots,-1$(resp. $+1$, $+2\cdots$, $+p,-p,\cdots$,-1) appear
in the same direction on all the vertices of $\Gamma_P$(resp.
$\Gamma_Q$). For example, the type B labels of Figure 5 are as in
Figure 6.

\begin{center}
\includegraphics[totalheight=4cm]{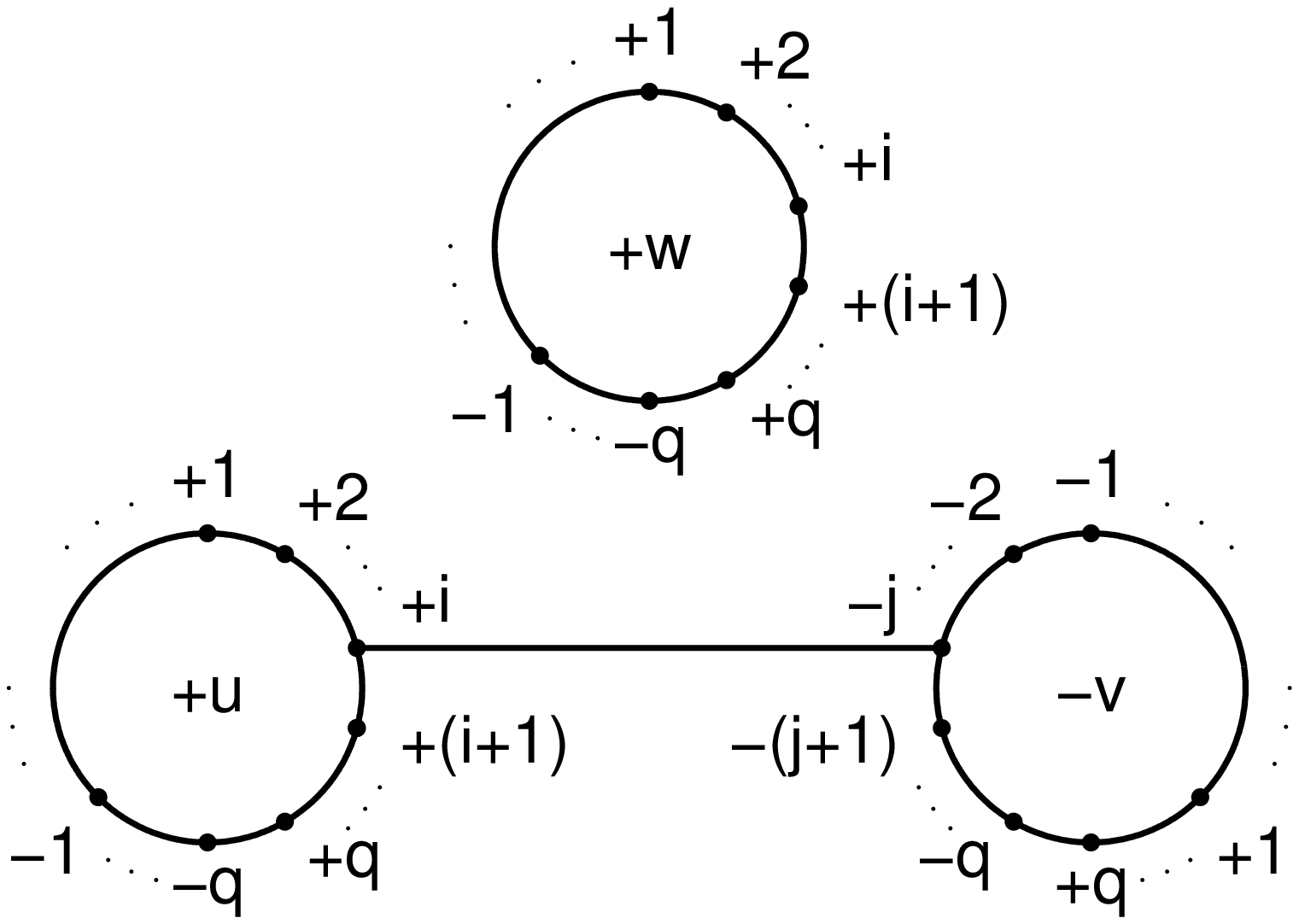}~~~~~~
\includegraphics[totalheight=4cm]{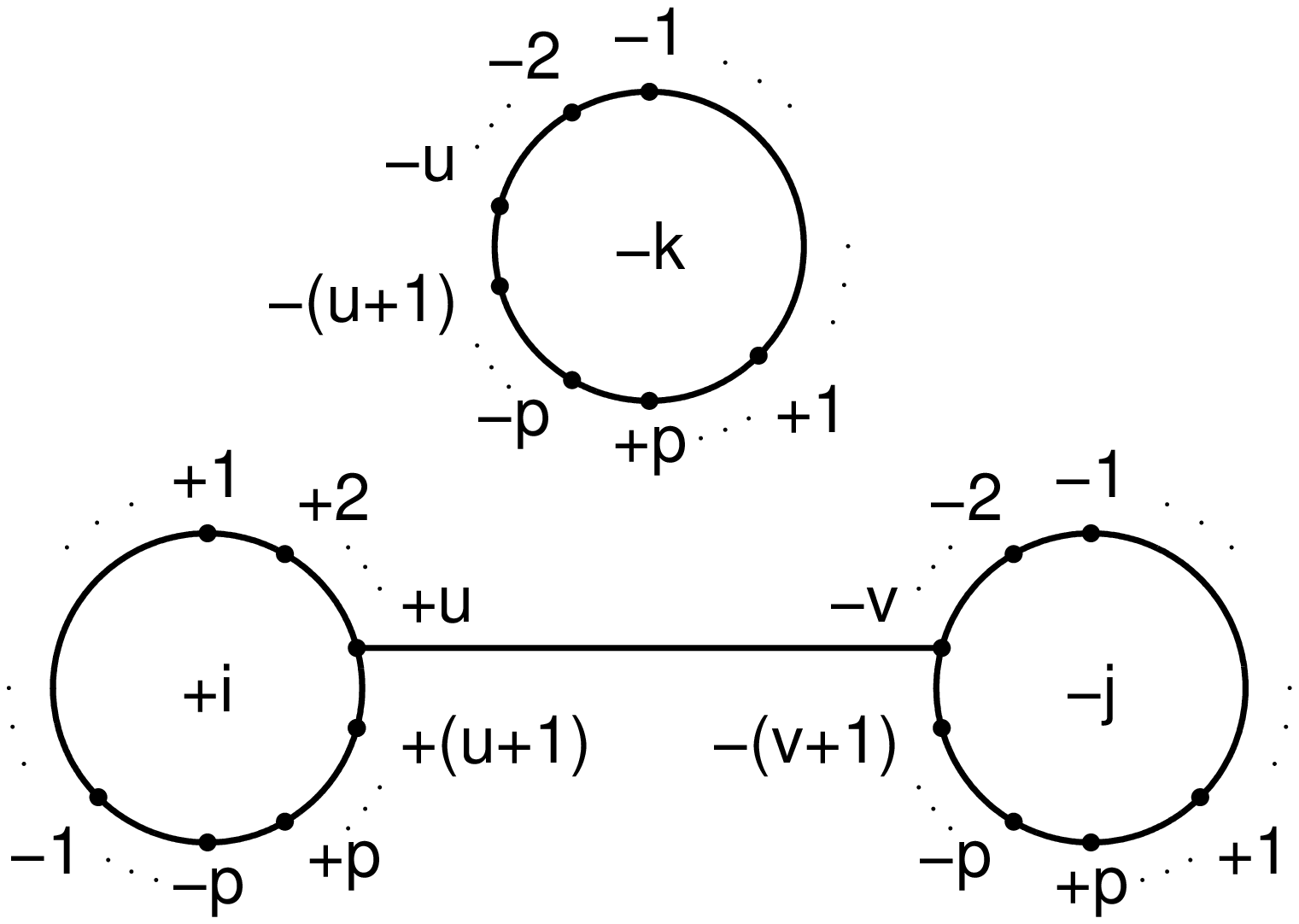}
\begin{center}
Figure 6
\end{center}
\end{center}

{\bf Lemma 3.2[ZQL].} \ (1) Let $e$ be an inner edge $e$ in
$\Gamma_P$ with its two endpoints $x$ labeled $(u,i)$ and $y$
labeled $(v,j)$, then $s(i)s(j)g(x)g(y)=-1$.

(2) \ Each inner edge in $\Gamma_P$ has its two endpoints with
different labels of type B. \hfill$\Box$\vskip 5mm

In the following arguments, the labels used for endpoints and edges
in $\Gamma_P$ and $\Gamma_Q$ are assumed to be type B.

{\bf Lemma 3.3.} \ Suppose $S=\{e_a \ | \ a=1,2,\cdots,m\}$ is a set
of $m$ parallel inner edges of $\Gamma_P$ with $e_a$ is labeled
$x_a-y_a$ as in Figure 7. If $x_a=y_b$ for some $1\leq a \leq m$ and
$1\leq b\leq m$, then $x_a=-y_a$ for each $1\leq a \leq m$.

{\bf Proof.} \ Without loss of generality, we assume that $a\leq b$.
By Lemma 3.2, $a<b$. By Remark 1, $x_{a+i}=y_{b-i}$. $b-a$ must be
odd, otherwise $x_{(a+b)/2}=y_{(a+b)/2}$ which contradicts Lemma
3.2. Hence $x_{(a+b-1)/2}=y_{(a+b+1)/2}$ and
$x_{(a+b+1)/2}=y_{(a+b-1)/2}$. This means that $e_{(a+b+1)/2}$ and
$e_{(a+b-1)/2}$ form a virtual Scharlemann cycle. By Lemma 2.8, it
is a virtual Scharlemann cycle rather than a Scharlemann cycle.
Hence $x_{(a+b-1)/2}=-y_{(a+b-1)/2}$, and $x_a=-y_a$ for each $1\leq
a \leq m$. \hfill$\Box$\vskip 5mm

{\bf Lemma 3.4.} \ Suppose $S=\{e_a \ | \ a=1,2,\cdots,m\}$ is a set
of $m$ parallel inner edges of $\Gamma_P$, and $e_a$ is labeled
$x_a-y_a$. If $m>q$, then $x_a=-y_a$ for each $1 \leq a \leq m$.

{\bf Proof.} \ If $m>q$, then there are some $a$ and $b$ such that
$x_a=y_b$. By Lemma 3.3, this lemma holds.\hfill$\Box$\vskip 5mm

\begin{center}
\includegraphics[totalheight=7cm]{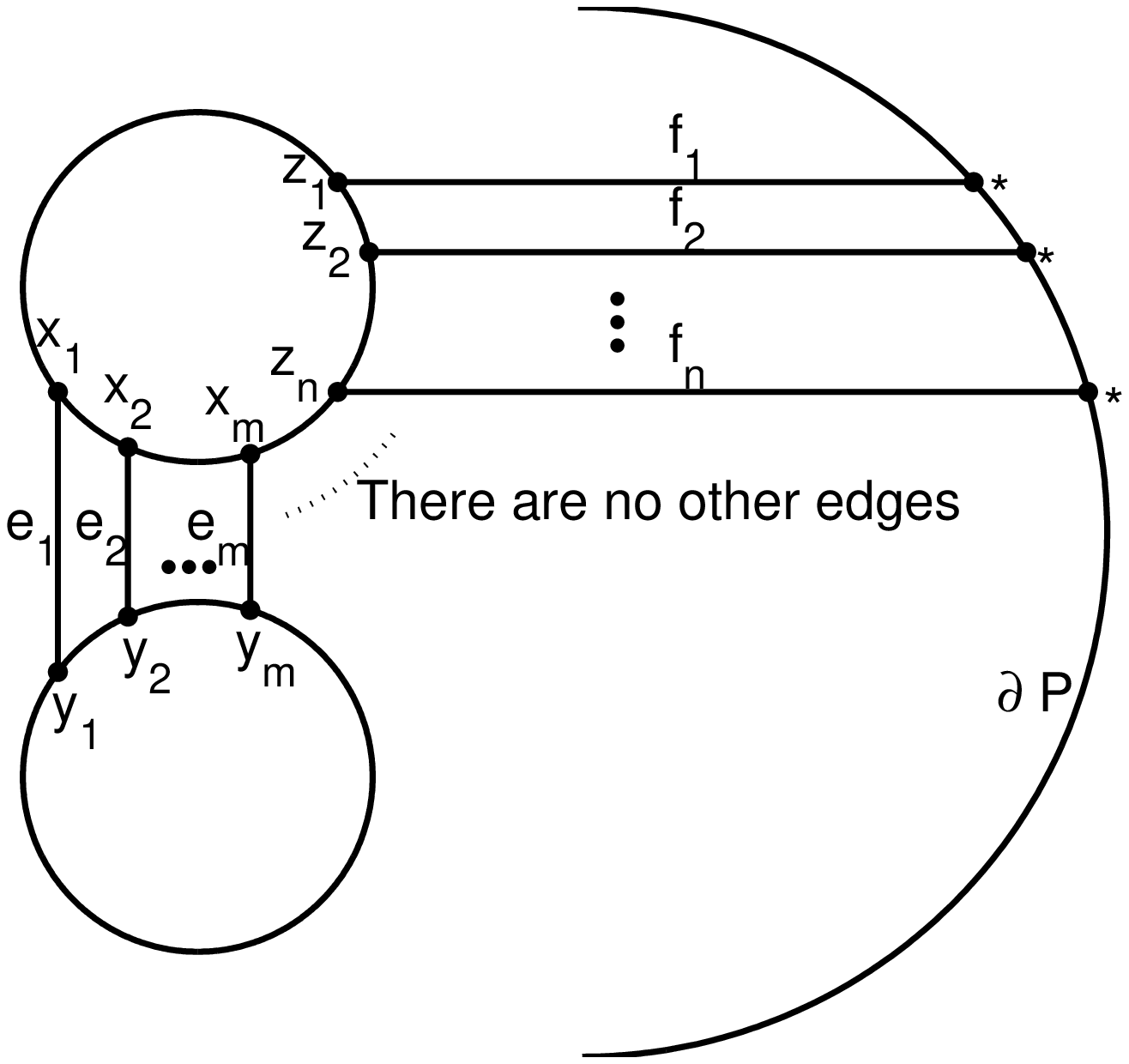}

Figure 7.
\end{center}

{\bf Lemma 3.5.} \ Let $S_1=\{e_a \ | \ a=1,2,\cdots,m\}$ be a set
of $m$ parallel inner edges, and $S_2=\{f_b\ | \  b=1,2\cdots,n\}$
be a set of $n$ parallel boundary edges of $\Gamma_P$ which is
adjacent to $S_1$. See Figure 7. If $n+m\geq2q$, then $m\leq q$.

{\bf Proof.} \ Suppose, otherwise, $m>q$. Let $e_a$ be labeled with
$x_a-y_a$ on $\Gamma_P$ as in Figure 7. By Lemma 3.4, $x_a=-y_a$.
Now, for each $1\leq a\leq m$, $e_a$ is a length one cycle in
$\Gamma_Q$ which is incident to $\partial_{|x_a|} S_\beta$ on
$\Gamma_Q$. Since $m+n\geq2q$, there is an edge in $S_2$ connecting
$\partial_x S_\beta$ to $\partial Q$ for $x\in
\{1,2,\cdots,q\}-\{|x_a| \ | \ a=1,2,\cdots,m\}$. This means that
the innermost one of the length one cycles $\{e_a\}$ is a trivial
loop in $\Gamma_Q$, a contradiction. \hfill$\Box$\vskip 5mm

For each $1\leq i \leq q$, let $B_P^{+i}$ be a subgraph of
$\Gamma_P$ consisting of all the vertices of $\Gamma_P$ and all the
edges $e$ such that one endpoint of $e$ is labeled with $+i$.

An $i$-triangle is a 3-sided disk face in $B_P^{+i}$. A boundary
$i$-triangle is an $i$-triangle such that one of its vertices is the
outer vertex of $\Gamma_P$.

The proof of Theorem 1 will be divided into two parts:

(1) \ $\Gamma_P$ has a boundary $i$-triangle.

(2) \ $\Gamma_P$ has no boundary $i$-triangle for each $1\leq i \leq
q$.

\begin{center}
\includegraphics[totalheight=4cm]{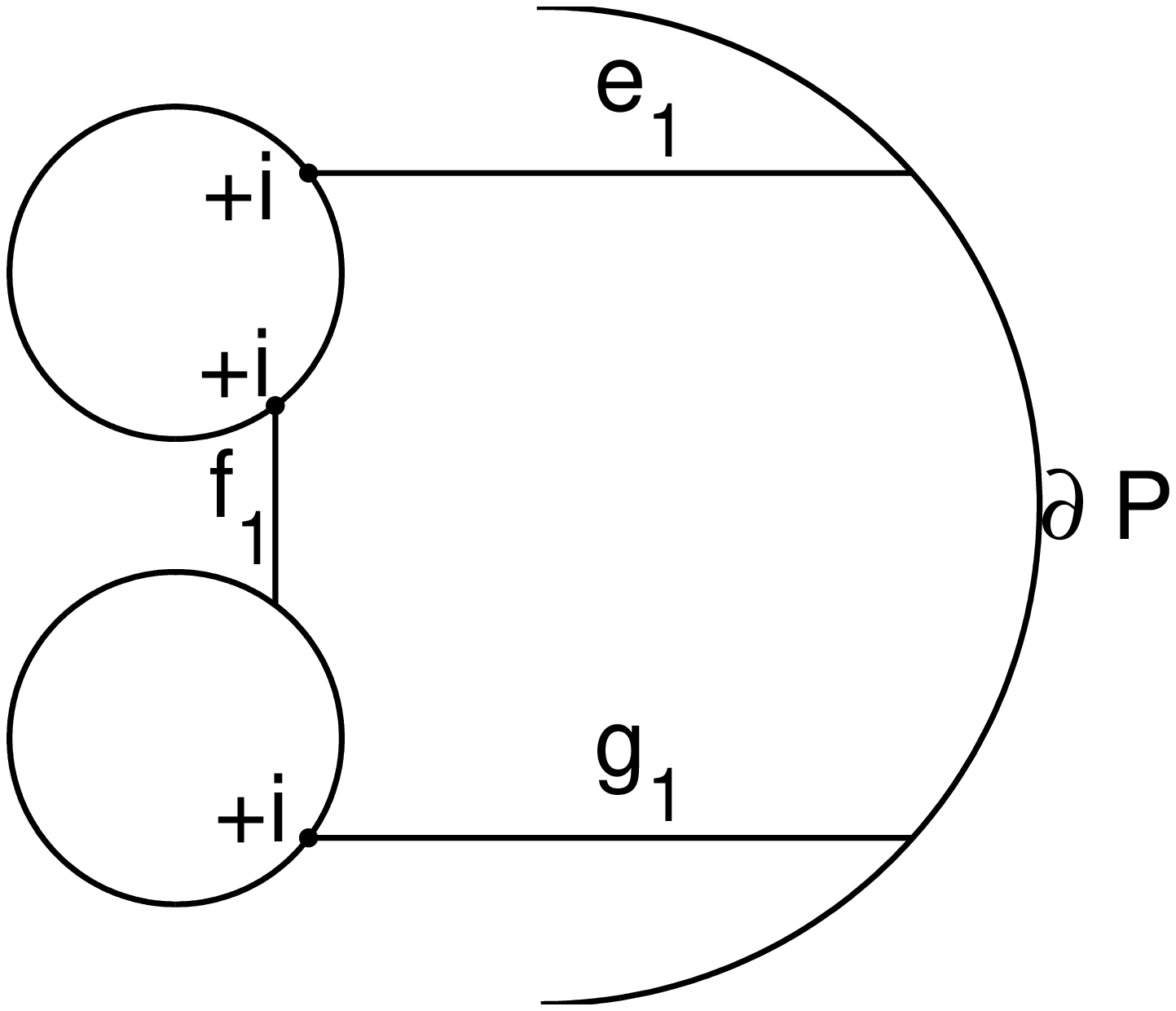}~~~~~
\includegraphics[totalheight=4cm]{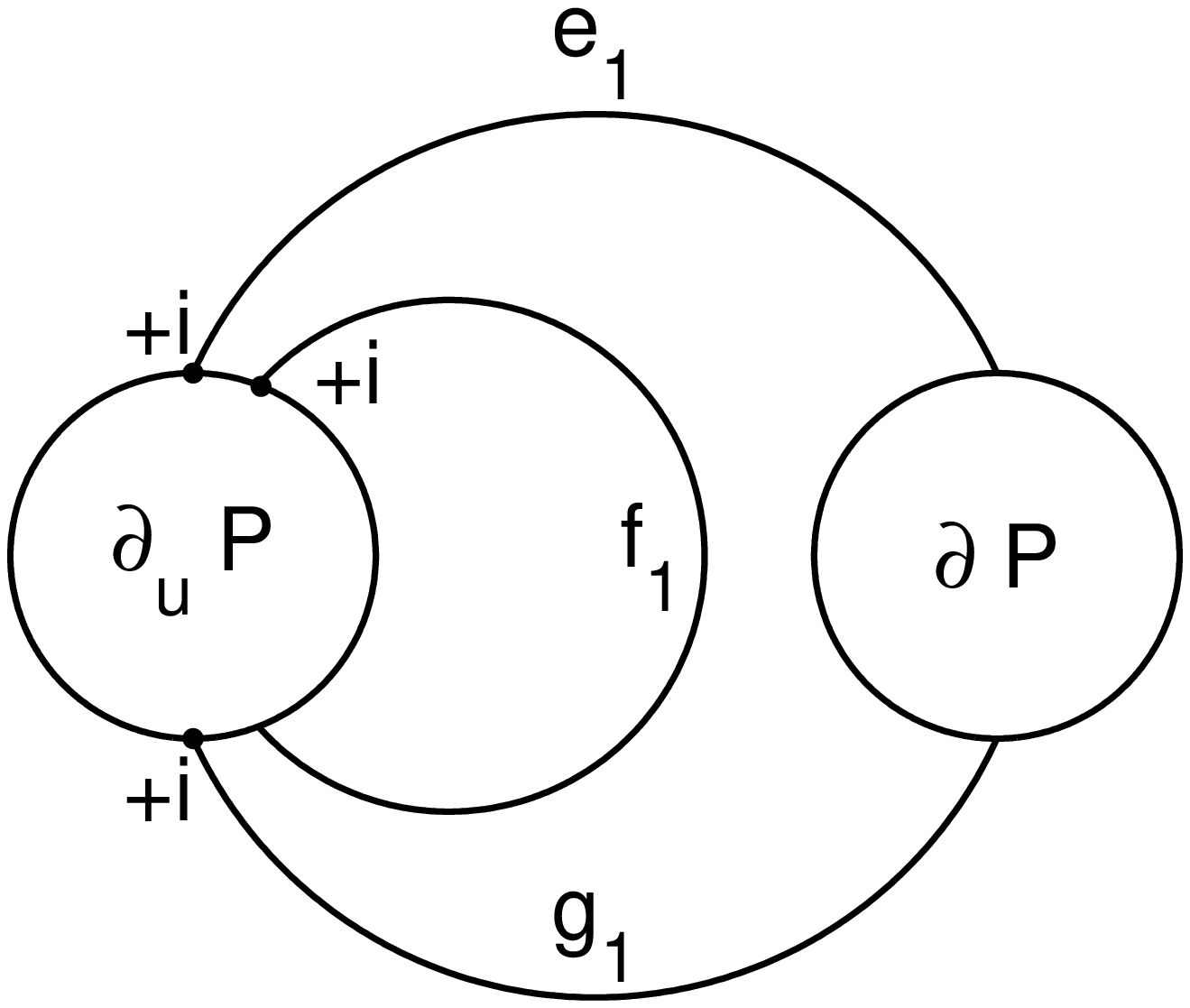}~~~~~
\includegraphics[totalheight=4cm]{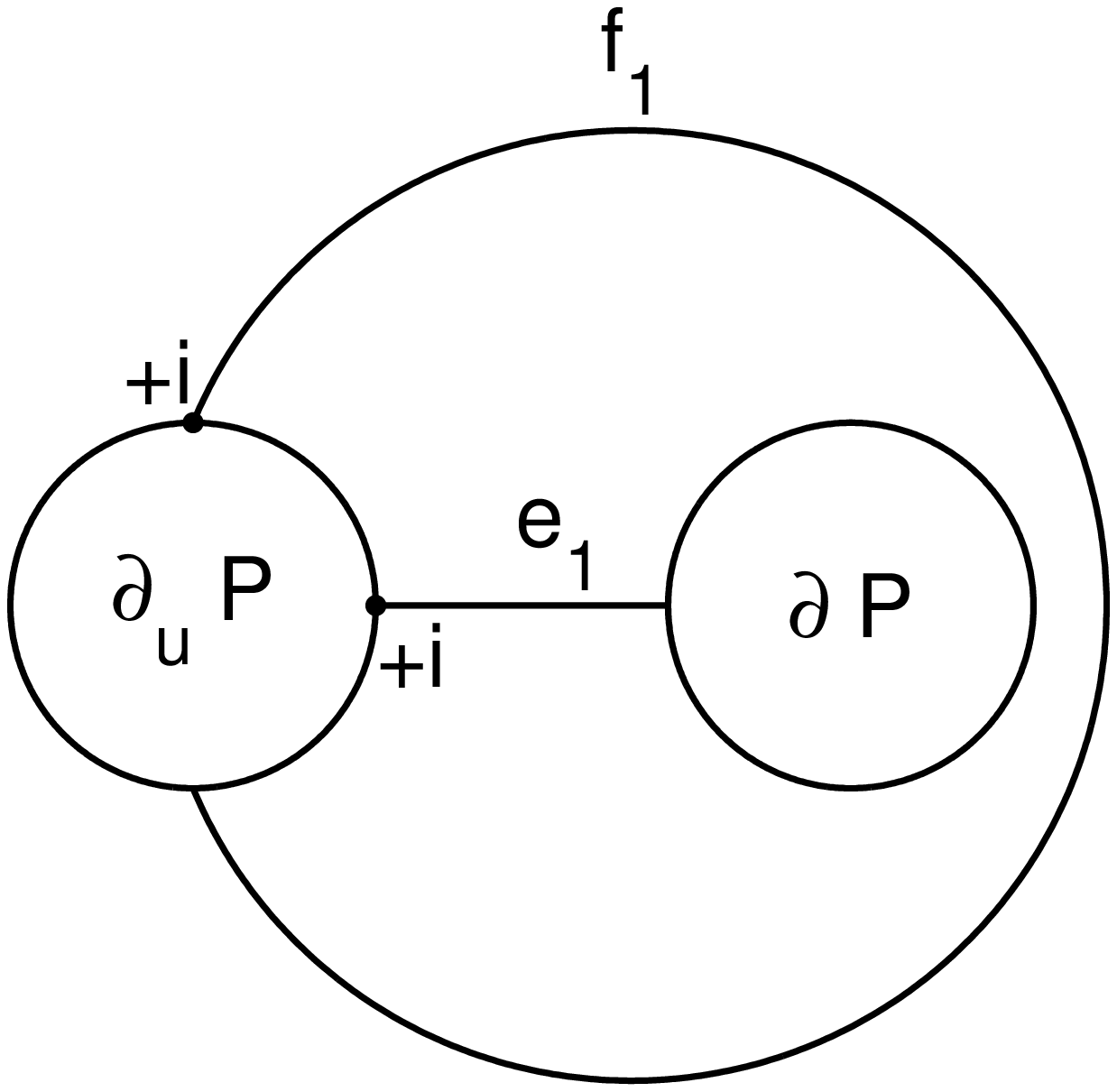}
(a)~~~~~~~~~~~~~~~~~~~~~~~~~~~~~~~(b)~~~~~~~~~~~~~~~~~~~~~~~~~~~~~~~(c)

Figure 8.
\end{center}

\begin{center}
\includegraphics[totalheight=8cm]{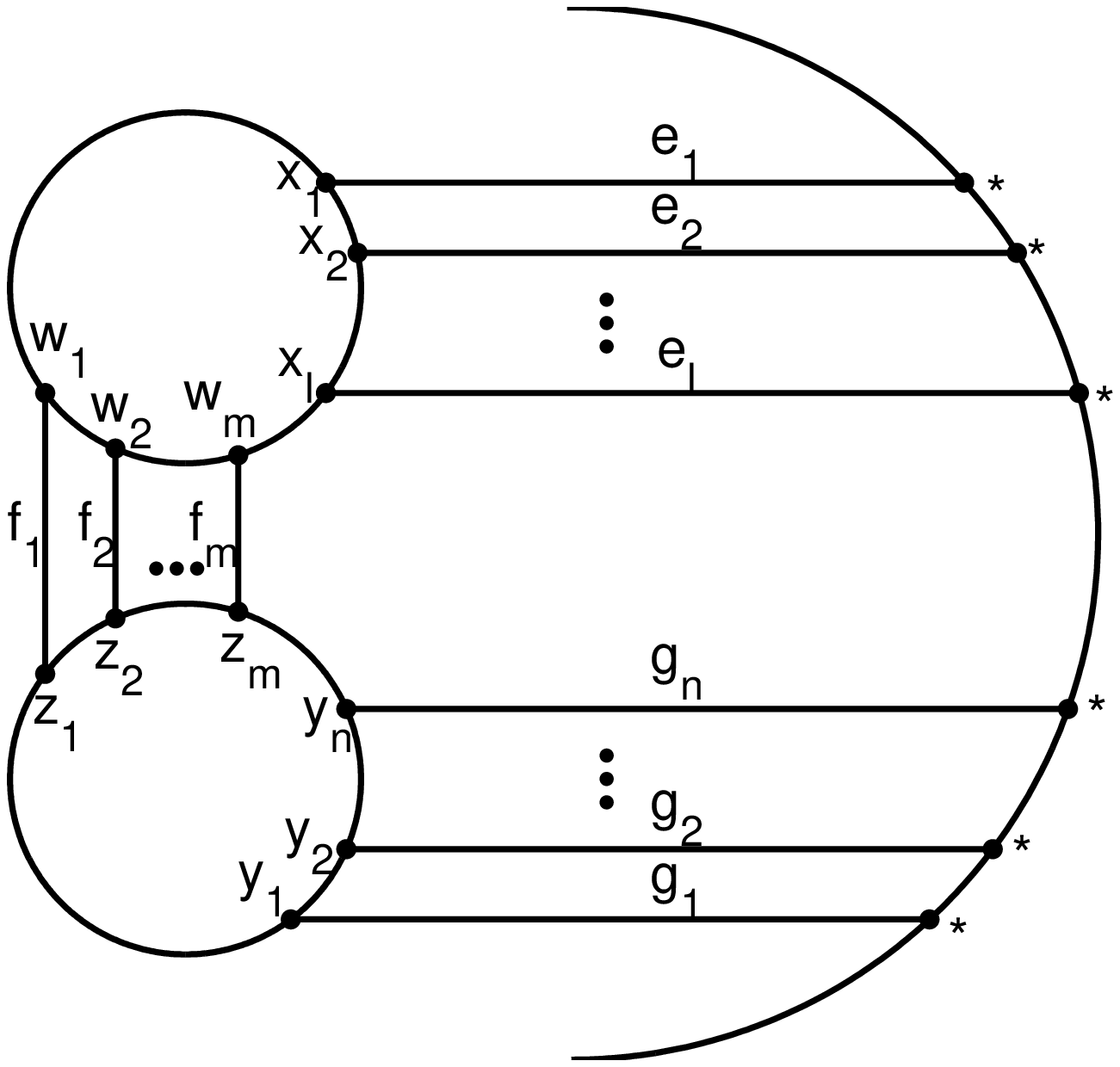}
\begin{center}
Figure 9.
\end{center}
\end{center}

\section {$\Gamma_P$ has a boundary $i$-triangle}

\ \ \ \ \ In this section, we assume $\Gamma_P$ has a boundary
$i$-triangle $\Delta$ for some $1\leq i\leq q$ as in Figure 8.

{\bf Lemma 4.1.} \ If $\Gamma_P$ contains a boundary $i$-triangle
$\Delta$ as in Figure 8, then each inner vertex of $\Gamma_Q$
belongs to a boundary edge.

{\bf Proof.} \  We first suppose $\Delta$ is as in one of Figure
8(a) and 8(b). There are three sets of parallel edges in $\Delta$.
We denote by $S_1=\{e_a \ | \ a=1,2,\cdots,l\}$ the set of the edges
parallel to $e_1$ in the $i$-triangle $\Delta$, $S_2=\{f_b \ | \
b=1,2,\cdots,m\}$ the set of the edges parallel to $f_1$ in
$\Delta$, and $S_3=\{g_c \ | \ c=1,2,\cdots,n\}$ the set of the
edges parallel to $g_1$ in $\Delta$. Furthermore, let $e_a$ be
labeled with $(u,x_a)-*$, $f_b$ be labeled with $(v,z_b)-(u,w_b)$
and $g_c$ be labeled $(v,y_c)-*$. Specially, $x_1=w_1=y_1=+i$. See
Figure 9.


Without loss of generality, we may assume that either $x_2> x_1$ or
$x_1=+q$ and $x_2=-q$. See Figure 9. Specially, if $x_1=+q$ and
$x_2=-q$, then we also assume that $x_2>x_1$.

Since $x_1=w_1=+i$, $l+m\geq 2q+1$. By Lemma 3.5, $m\leq q$. If
$w_m>0$, then the labels $-q,-(q-1),\cdots,-1$ appear in $\{x_a \ |
\ a=1,2,\cdots,l\}$. This means that, for each vertex $v$ of
$\Gamma_Q$, there is an edge in $S_1$ connecting $v$ to the outer
vertex.

Now assume that $w_m<0$. Since $w_1=+i$ and $x_2>x_1$, $w_2<w_1$.
Hence $+1\in \{w_b \ | \ b=1,2,\cdots,m\}$. By Remark 1 in Section
3, the labels $+1$, $+2\cdots$, $+q,-q,\cdots,-1$ appear in the same
direction on all the vertices of $\Gamma_P$. Hence $y_2<y_1=+i$.
This means that $+1\in \{z_b,y_c \ \ | \ 1\leq b\leq m; 1\leq c\leq
n\}$.

Case 1. $+1\in \{z_b \ | \ 1\leq b\leq m\}$.

By Lemma 3.3, $w_b=-z_b$. Since $l+m\geq 2q+1$, by the proof of
Lemma 3.7, it is impossible.

Case 2. $+1\in\{y_c \ | \ 1\leq c\leq n\}$.

Since $w_m<0$ and $x_2>x_1=+i$, we have
$\{+i,+(i+1),\cdots,+q\}\subset \{x_a \ | \ a=1,2,\cdots,l$\}. Since
$y_2<y_1=+i$, $\{+1,+2,\cdots,+i\}\subset\{y_c \ | \ c=1,2,\cdots, n
\}$. Hence each inner vertex of $\Gamma_Q$ belongs to a boundary
edge.

Suppose now that $\Delta$ is as in Figure 8(c). In this case, we can
also find three sets of parallel edges in $\Delta$ as in Figure 10.
Specially, let $g_1=e_1$. By the same argument as above, the lemma
holds.\hfill$\Box$\vskip 5mm
\begin{center}
\includegraphics[totalheight=7cm]{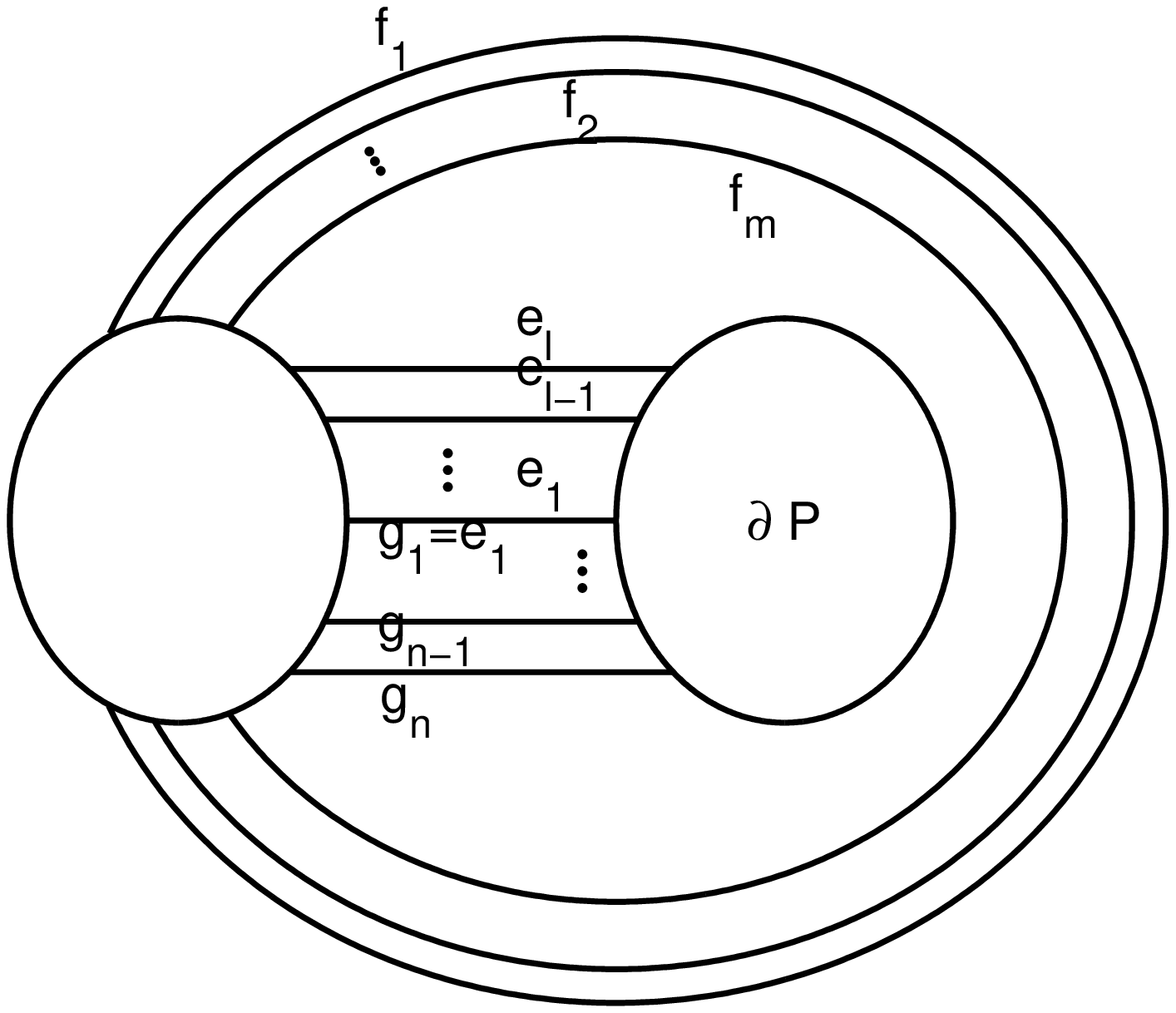}
\begin{center}
Figure 10.
\end{center}
\end{center}

We denote by $\bar\Gamma_Q$ the reduced graph of $\Gamma_Q$.

{\bf Lemma 4.2.} \ $\bar{\Gamma}_Q$ has a vertex of valence at most
three which belongs to a single boundary edge.

{\bf Proof.} \ The lemma follows from Lemma 4.1 and Lemma 2.6.5 in
[CGLS].\hfill$\Box$\vskip 5mm

{\bf Proposition 4.3.} Theorem 1 is true for the case: $\Gamma_P$
contains a boundary $i$-triangle.

{\bf Proof.} \ By Lemma 4.2, $\bar\Gamma_Q$ contains a vertex of
valence at most three, say $v$. Since $\Delta(\alpha,\beta)\geq4$,
by lemma 2.9, $v$ is of valence three. Hence there are three sets of
parallel edges incident to $v$, say $\{e_a \ | \ a=1,2,\cdots,l\}$,
$\{f_b \ | \ b=1,2,\cdots,m\}$ and $\{g_c \ | \ c=1,2,\cdots,n\}$.
Without loss of generality, we may assume that $\{e_a\}$ are
boundary edges while $\{f_b\}$ and $\{g_c\}$ are inner edges. Hence
$l+m+n\geq4p$. By Lemma 2.9, $l<2p$. Hence one of $m$ and $n$, say
$m>p$. Also by Lemma 2.9, $n<2p$. Hence $l+m>2p$. It is a
contradiction to Lemma 3.5. \hfill$\Box$\vskip 5mm

\section {$\Gamma_P$ has no boundary i-triangle}

\ \ \ \ \ In this section, we assume that $\Gamma_P$ contains no
boundary $i-$triangle for each $1\leq i \leq q$.

{\bf Lemma 5.1} \ $\Gamma_P$ has an edge labeled with $(+i)-(-i)$
for each $1\leq i\leq q$.

{\bf Proof.} \ Since $\Delta(\alpha,\beta)\geq 4$, by Lemma 3.2(2),
$B_P^{+i}$ has at least $2p$ edges.

{\bf Claim.} \ $B_P^{+i}$ has at least one 2-sided or 3-sided disk
face.

{\bf Proof. } \ Denoted by $V$, $E$ and $F$ the number of vertices,
edges and disk faces in $B_P^{+i}$. Then $V=p+1$, and $V\leq E/2+1$.
(In this case, we take $B_P^{+i}$ as a graph in a 2-sphere.)

Suppose, otherwise, that $B_P^{+i}$ contains no 2-sided and 3-sided
disk faces. Then $2E\geq 4F$. Hence $V-E+F<E/2+1-E+E/2=1<2$, a
contradiction.  \hfill$\Box$(Claim )\vskip 5mm

Now if $B_P^{+i}$ contains a 2-sided or 3-sided disk face, say $C$,
then, by Lemma 2.9 and Proposition 4.3, all the edges in $C$ are
inner.

Case 1. \ $B_P^{+i}$ contains a 2-sided disk face.

Now this 2-sided disk face offers $n$ adjacent parallel edges
$e_1,e_2,\cdots e_n$ of $\Gamma_P$ with $e_a$ labeled with
$x_a-y_a$, such that $x_1=+i$ and one of $x_n$ and $y_n$ is also
$+i$. By Lemma 2.9, $y_n=+i$. By Lemma 3.3, $x_a=-y_a$ for each
$1\leq a\leq n$. Hence $e_1$ is labeled with $(+i)-(-i)$.

Case 2. \ There is a 3-side disk face $\Delta$ in $B_P^{+i}$.

Now $\Delta$ is as in one of Figures 5.1(a), (b) and (c). If
$\Delta$ is as in Figure 11(c), then $\Gamma_P$ contains $4q$
parallel edges, contradicting Lemma 2.9.


\begin{center}
\includegraphics[totalheight=4cm]{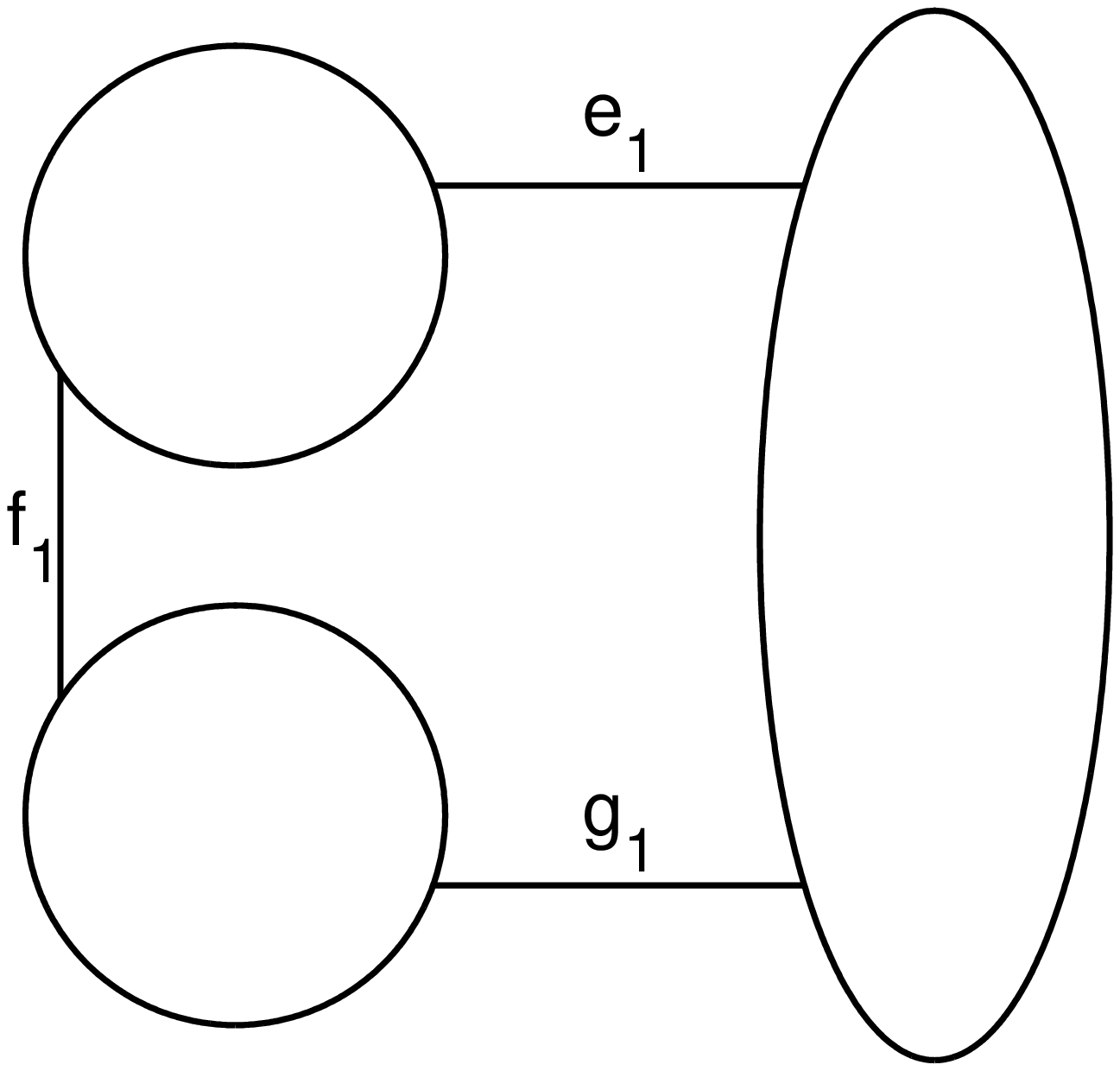}~~~~~
\includegraphics[totalheight=4cm]{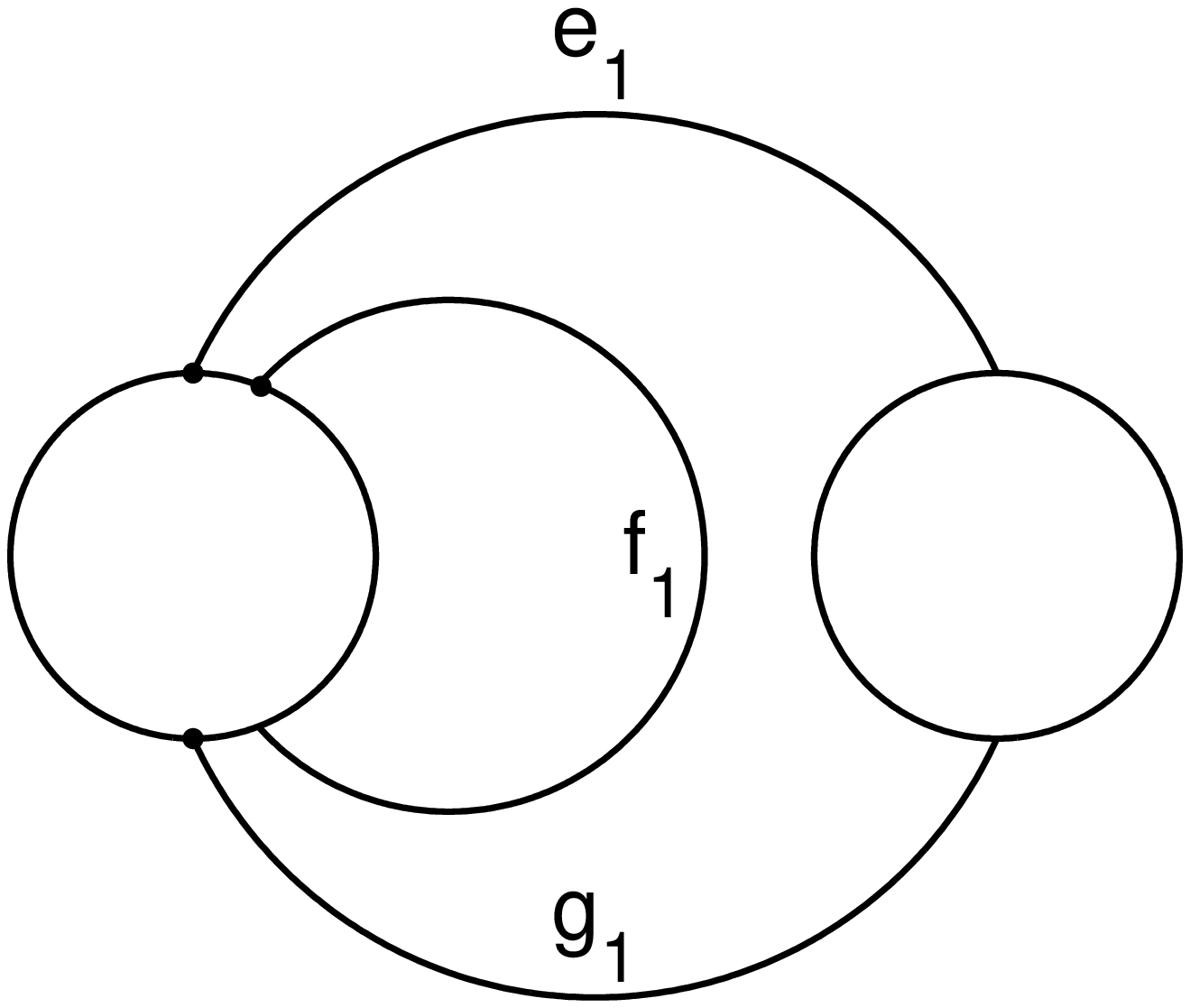}~~~~~
\includegraphics[totalheight=4cm]{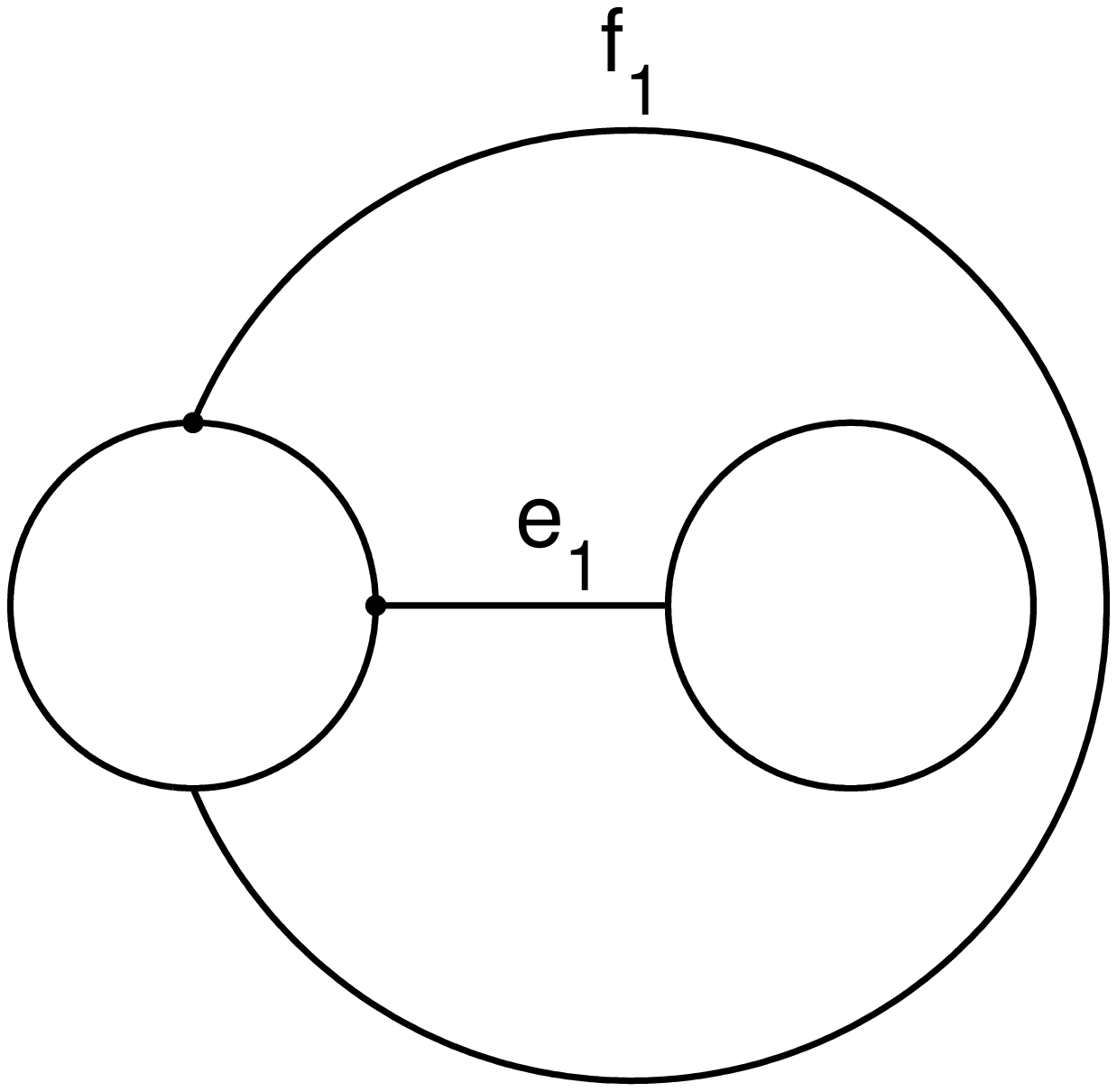}
(a)~~~~~~~~~~~~~~~~~~~~~~~~~~~~~~~(b)~~~~~~~~~~~~~~~~~~~~~~~~~~~~~~~(c)

Figure 11.
\end{center}

Suppose now that $\Delta$ is as in Figure 11(a). We denote by $e_1$,
$f_1$ and $g_1$ the three boundary edges of $\Delta$. We denote by
$S_1=\{e_a \ | \ a=1,2,\cdots,l\}$ the set of the edges parallel to
$e_1$ in the $i$-triangle $\Delta$, $S_2=\{f_b \ | \
b=1,2,\cdots,m\}$ the set of the edges parallel to $f_1$ in
$\Delta$, and $S_3=\{g_c | \ c=1,2,\cdots,n\}$ the set of the edges
parallel to $g_1$ in $\Delta$. See Figure 12. Let $e_a$ be labeled
with $x_a-y_a$, $f_b$ be labeled with $z_b-w_b$, and $g_c$ be
labeled with $s_c-t_c$.

\begin{center}
\includegraphics[totalheight=8cm]{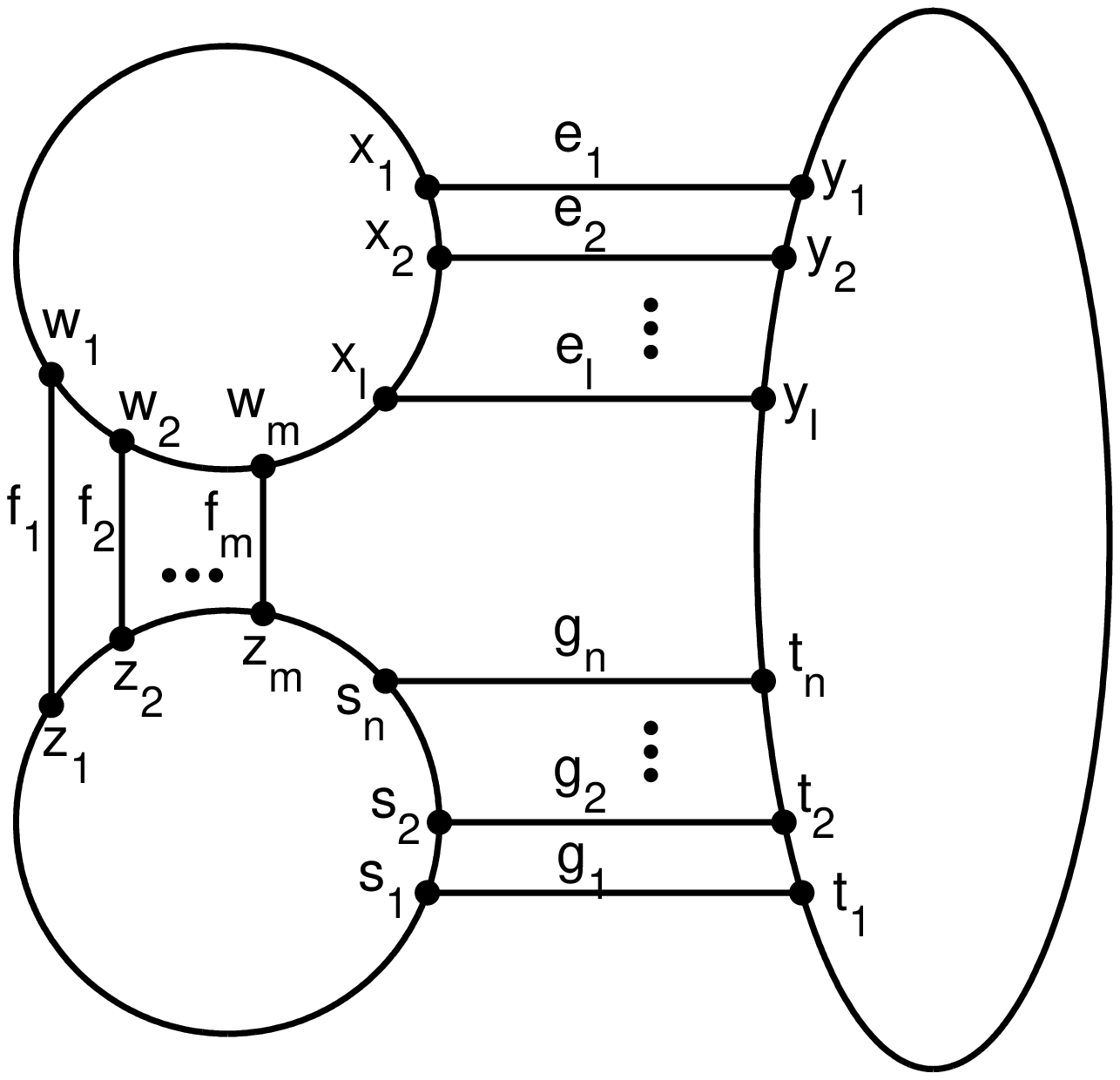}
\begin{center}
Figure 12.
\end{center}
\end{center}

\begin{center}
\includegraphics[totalheight=4cm]{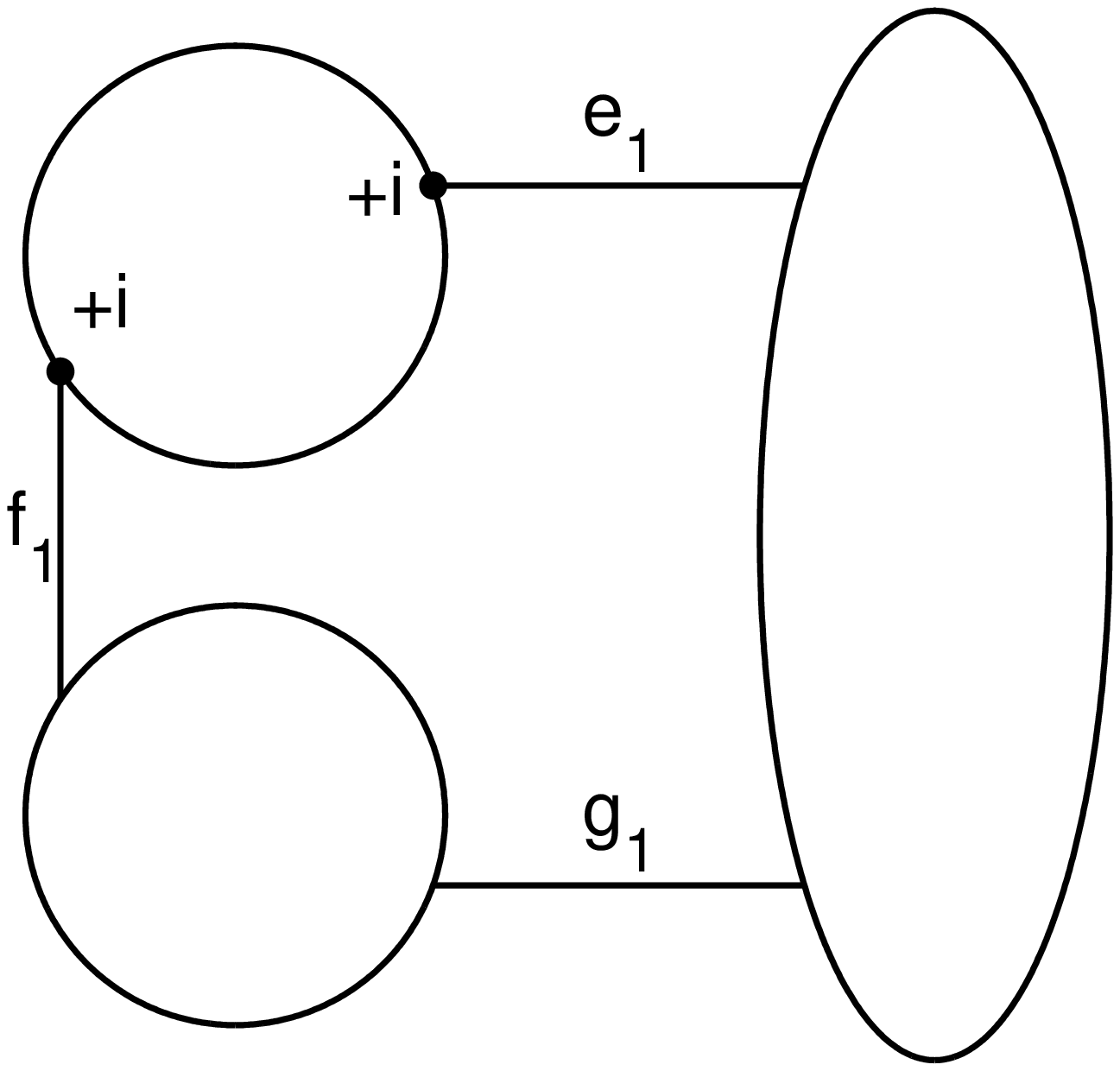}~~~~~
\includegraphics[totalheight=4cm]{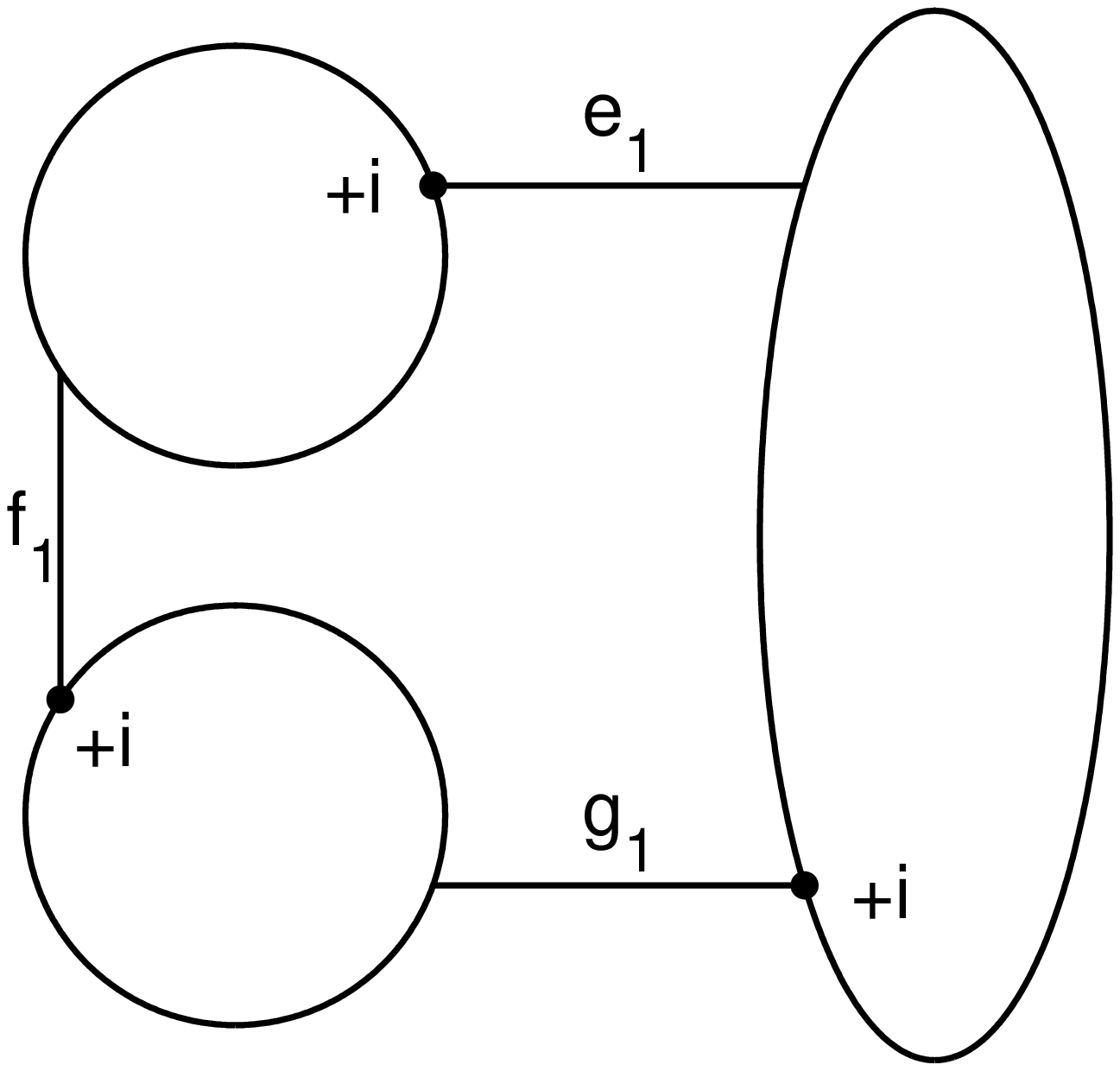}~~~~~

(a)~~~~~~~~~~~~~~~~~~~~~~~~~~~~~~~(b)

Figure 13.
\end{center}

Case 2.1. \ $x_1=w_1=+i$ as in Figure 13(a).



Now $l+m>2q$, this indicates that one of $l$ or $m$, say $l>q$. By
Lemma 3.4, $x_a=-y_a$ for each $1\leq a\leq l$. Hence $e_1$ is the
edge labeled with $(+i)-(-i)$.

Case 2.2. \ $x_1=z_1=t_1=+i$ as in Figure 13(b).

We first suppose $l=m=n$. By Remark 1 in Section 3, $x_l=z_m=t_n$.
Hence $y_l=w_m=s_n$. This means that $e_l$, $f_m$ and $g_n$ form a
virtual Scharlemann cycle. By Lemma 2.8, $x_l=-y_l=\pm1$ or
$x_l=-y_l=\pm q$. Hence $x_a=-y_a$ for each $1 \leq a\leq n$, and
$e_1$ is labeled with $(+i)-(-i)$.

Now we suppose $l<m$. By Remark 1 in Section 3, $z_a=x_a$ for each
$1\leq a \leq l$, and $z_{l+1}=w_m$. By Lemma 3.3, $z_b=-w_b$ for
each $1\leq b \leq m$. Then $w_1=-z_1=-i$. Hence $f_1$ is labeled
with $(+i)-(-i)$. \hfill$\Box$\vskip 5mm

{\bf The proof of Theorem 1.} \ By Lemma 5.1 and Proposition 4.3,
for each $1\leq i\leq q$, there is an edges labeled with
$(+i)-(-i)$. This means that for each vertex $\partial_i S_\beta$ of
$\Gamma_Q$, there is an edge with its two endpoints incident to
$\partial_i S_\beta$. Hence $\Gamma_Q$ contains a 1-sided disk face,
a contradiction. By Lemma 2.3, $\alpha=\beta$. Hence Theorem 1
holds.\hfill$\Box$\vskip 5mm

{\bf The proof of Theorem 2.} \ Under the assumptions of Theorem 2,
all the arguments except Lemma 2.3 are true. Hence Theorem 2
holds.\hfill$\Box$
\noindent Yannan Li \\
Department of Applied Mathematics, \\Dalian University of Technology,\\
Dalian, China, 116022\\ Email: yn\_lee79@yahoo.com.cn

\vskip 5mm
\noindent Ruifeng Qiu\\
Department of Applied Mathematics, \\Dalian University of
Technology,\\
Dalian, China, 116022  \\Email: qiurf@dlut.edu.cn

\vskip 5mm

\noindent Mingxing Zhang \\
Department of Applied Mathematics, \\Dalian University of
Technology,
\\Dalian, China, 116022 \\Email: star.michael@263.net

\end{document}